\documentclass[dvipsnames]{trbunofficial}

\usepackage{graphicx}
\usepackage{tikz}
\usepackage{bbm}
\usepackage[ruled, linesnumbered]{algorithm2e}
\SetAlCapNameFnt{\scriptsize}
\SetAlCapFnt{\footnotesize}
\usepackage[fleqn]{amsmath}
\usepackage{algpseudocode}
\usepackage{float}
\usepackage{setspace}
\usepackage[T1]{fontenc}

\usepackage[colorlinks=true,linkcolor=blue,citecolor=blue]{hyperref}
\usepackage{amsthm,amsmath,amssymb}
\usepackage{mathtools}
\usepackage{cleveref}
\usepackage[font=small,labelfont=bf]{caption}
\usepackage{subcaption}

\usepackage[utf8]{inputenc}
\usepackage[english]{babel}

\usepackage{amsthm}
\usepackage[pagewise]{lineno}

\let\oldequation\equation
\let\oldendequation\endequation

\renewenvironment{equation}
  {\linenomathNonumbers\oldequation}
  {\oldendequation\endlinenomath}

\newcount\Comments  
\newcommand{\kibitz}[2]{\ifnum\Comments=0\textcolor{#1}{#2}\fi}

\newtheorem{definition}{Definition}[section]

\newtheorem{property}{Property}[section]

\newtheorem*{theorem*}{Theorem}

\newtheorem{lemma}{Lemma}[section]
\newtheorem*{conjecture*}{Conjecture}

\numberwithin{equation}{section}
\theoremstyle{remark}
\newtheorem{example}{Example}[section]

\newcommand{\inv}{\ensuremath{^{-1}}}
\newcommand{\transpose}{\ensuremath{^\top}}
\newcommand{\optimal}{\ensuremath{^*}}
\newcommand{\nmpc}{\ensuremath{{n_{e}T}}}
\newcommand{\mmpc}{\ensuremath{{n_{e}}}}
\newcommand{\by}{\ensuremath{\times}}
\newcommand{\qmpc}{\ensuremath{q}}

\newcommand{\R}{\ensuremath{\mathbb{R}}}

\newcommand{\stateset}{\ensuremath{\mathcal{X}}}
\newcommand{\edgeset}{\ensuremath{\mathcal{E}}}
\newcommand{\network}{\ensuremath{\mathcal{G}}}
\newcommand{\vertexset}{\ensuremath{\mathcal{V}}}
\newcommand{\edgeflowset}{\ensuremath{\mathbb{R}_{\geq 0}}}
\newcommand{\costset}{\ensuremath{\mathbb{R}_{\geq 0}}}

\newcommand{\timesym}{\ensuremath{t}}

\newcommand{\edgeflow}{\ensuremath{f_e}}
\newcommand{\edgeelem}{\ensuremath{e}}
\newcommand{\numedges}{\ensuremath{n_e}}
\newcommand{\flow}{\ensuremath{f}}
\newcommand{\latfunc}{\ensuremath{l_e}}
\newcommand{\flowset}{\ensuremath{\Delta}}

\newcommand{\flowstoch}{\ensuremath{f}}
\newcommand{\otheredgeflow}{\ensuremath{f_{e'}}}
\newcommand{\otheredgelatfunc}{\ensuremath{l_{e'}}}
\newcommand{\otheredge}{\ensuremath{e'}}
\newcommand{\origin}{\ensuremath{o}}
\newcommand{\destination}{\ensuremath{d}}

\newcommand{\loss}{\ensuremath{l}}
\newcommand{\flowvectorset}{\ensuremath{\edgeflowset^{\numedges}}}
\newcommand{\flowiter}{\ensuremath{\flowstoch_t}}
\newcommand{\flowiterplus}{\ensuremath{\flowstoch_{t+1}}}
\newcommand{\lossiter}{\ensuremath{\loss_t}}
\newcommand{\lossvectorset}{\ensuremath{\R^{\numedges}}_{\geq 0}}
\newcommand{\flowdyn}{\ensuremath{h}}
\newcommand{\flowcostfunciter}{\ensuremath{c_{2}}}

\newcommand{\state}{\ensuremath{\mathbf{x}}}
\newcommand{\control}{\ensuremath{\mathbf{u}}}

\newcommand{\F}{\ensuremath{F}}
\newcommand{\Y}{\ensuremath{Y}}
\newcommand{\stateiter}{\ensuremath{\state_t}}
\newcommand{\controliter}{\ensuremath{\control_t}}
\newcommand{\stateiterplus}{\ensuremath{\state_{t+1}}}
\newcommand{\ctrldyn}{\ensuremath{g}}
\newcommand{\costfunciter}{\ensuremath{c_{1}}}
\newcommand{\costfunciterT}{\ensuremath{c_{1,T}}}
\newcommand{\initstate}{\ensuremath{\state_0}}

\newcommand{\nU}{\ensuremath{n_u}}
\newcommand{\nx}{\ensuremath{n_x}}
\newcommand{\statevectorset}{\ensuremath{\R^{\nx}}}
\newcommand{\controlvectorset}{\ensuremath{\R^{\nU}}}
\newcommand{\edgestate}{\ensuremath{\state_{e}}}
\newcommand{\otheredgestate}{\ensuremath{\state_{e'}}}
\newcommand{\otheredgestatelatfunc}{\ensuremath{l_{e'}}}

\newcommand{\CR}{\ensuremath{\mathcal{R}}}
\newcommand{\mpcctrl}{\ensuremath{\mathcal{F}}}
\newcommand{\mpcctrladd}{\ensuremath{\mathcal{G}}}
\newcommand{\CRCheckH}{\ensuremath{\mathcal{H}}}
\newcommand{\CRCheckK}{\ensuremath{\mathcal{K}}}
\newcommand{\condenseqpobj}{\ensuremath{V}}
\newcommand{\activeset}{\ensuremath{\mathcal{A}}}

\newcommand{\regionelem}{\ensuremath{r}}

\AuthorHeaders{Gibson, You, and Bayen}
\title{Parallel Network Flow Allocation in Repeated Routing Games via LQR Optimal Control}

\author{%
  \textbf{Marsalis Gibson*}\\
  University of California, Berkeley\\
  652 Sutardja Dai Hall, Berkeley, CA 94720-1710\\
  email: mtgibson@berkeley.edu\\
  \hfill\break
  \textbf{Yiling You*}\\
  University of California, Berkeley\\
  652 Sutardja Dai Hall, Berkeley, CA 94720-1710\\
  email: yiling.you@berkeley.edu \\
  \hfill\break
  \textbf{Alexandre M. Bayen}\\
  University of California, Berkeley\\
  Institute of Transportation Studies\\
  109 McLaughlin Hall, Berkeley CA 94720-1720\\
  email: bayen@berkeley.edu \\
  \hfill\break
}

\usepackage{todonotes}

\begin{document}

\maketitle

\section{Abstract}

In this article, we study the repeated routing game problem on a parallel network with affine latency functions on each edge. We cast the game setup in a LQR control theoretic framework, leveraging the Rosenthal potential formulation. We use control techniques to analyze the convergence of the game dynamics with specific cases that lend themselves to optimal control. We design proper dynamics parameters so that the conservation of flow is guaranteed. We provide an algorithmic solution for the general optimal control setup using a multiparametric quadratic programming approach (explicit MPC). Finally we illustrate with numerics the impact of varying system parameters on the solutions.

\hfill\break%
\noindent\textit{Keywords}: Repeated routing game, Algorithmic Game Theory, Nash Equilibrium, Linear-quadratic Regulator, Multiparametric quadratic programming, Explicit MPC 

\newpage

\section{Introduction}
The routing problem -- a problem in which trip-makers route traffic between an origin and destination -- is a central and essential component in the study of urban transportation planning. Emerging as an enabler to this problem over the last decade, navigational apps have used solution from the routing problem to provide shortest paths from origin to destination, or identify the paths with the shortest travel time \cite{Cabannes2018TheIO}. Due to ubiquitous adoption of navigational apps by motorists (commuters, travelers), TNCs (Lyft/Uber), delivery fleets \cite{hwang2015urban}, and many others, these apps are playing even bigger roles in transportation and the significance of their roles are growing \cite{alonso2017demand}. With increased penetration routing app users, navigational apps are starting to "act" as flow allocation mechanisms -- having positive effects on some cities travels and negative on others when users use local neighborhood routes to avoid traffic \cite{thai2016negative}. Thus, it is becoming ever for important to understand the routing game problem, its solutions, and how flow allocation mechanisms (and traffic assignment) are being used in transportation. 

The traffic assignment problem is a model framework for aggregated flows of motorists routed through a network from their origins to their destinations. To address the traffic assignment problem, this process has initially been modeled in a static way using notions such as static traffic assignment, for which specific solutions exist like User equilibrium and system optimal \cite{boyles2010traffic} \cite{patriksson2015traffic}. In this static model, we assume the network being studied is close to "steady state", and that the current network capacity is the capacity over the entire analysis period. However, in hopes to move away from this assumption and better model traffic in urban areas, researchers have moved toward dynamic traffic assignment models, which aim to
model time-varying networks and demand interactions and changes \cite{chiu2011dynamic} \cite{shafiei2018calibration} \cite{mahmassani2001dynamic} \cite{boyles2010traffic}. Additionally, within the field of traffic control, researchers have used simulation based approaches \cite{osorio2010simulation} and some are starting to use big data to further analyze travel behavior  \cite{chen2016promises}. 

The traffic assignment problem has been modeled using game theory among many other techniques. Specifically, it has been posed in a game theoretic framework called routing game. A routing game is a game in which the players (non-cooperatively) route traffic between an origin and a destination of a proposed network, where the network models a system of roads. After each player chooses routing strategies, each player will incur a cost on the route they choose \cite{Nisam2007}. One specific kind of game is called "the one-shot" game. The one-shot routing game is a static game\footnote{A static game is played only once with one iteration.}, where the players make their decision simultaneously and therefore no time history is involved. Within the one-shot game framework, multiple equilibrium models and the corresponding equilibria have played an important role in understanding the inefficiencies of networks (e.g. the price of anarchy \cite{10.1145/506147.506153}), for example the user equilibrium or Nash equilibrium \cite{patriksson}. Properties and the characterizations of the equilibrium solutions are also well-studied, for example, the Nash equilibria are known to be the solution set of the convex Rosenthal potential function \cite{Rosenthal1973ACO}. 

Within routing games, one can also study a dynamic game, where the environment is time varying. The analysis of such games is more difficult from a theoretic standpoint: the notion of dynamic equilibria has to be introduced, and extra ingredients like information acquisition and time ordering of the execution of the actions need to be specified \cite{doi:10.1137/1.9781611971132}. Therefore, moving towards more dynamical games, we consider repeated games, i.e. a static game iterated over time, as a model of the evolution of the traffic situation. The iteration here can be naturally interpreted as time discretization (e.g. sometimes where traffic conditions are steady over time intervals). In the case of the present work, we follow previous approaches in which the repeated nature of the game encompasses the process of a routing entity (Google, Waze, etc) to learn from day to day from their previous performances. Oftentimes (and it is the case in our formulation), the repeated version is essentially the sequential appreciation of a static game. A natural question to study is the behavior of the dynamics, when the players repeat a static game using past information. Dynamics in a repeated routing game might not converge, and even if it does, convergence might not be quantifiable (e.g. convergence almost everywhere, convergence in the sense of Cesaro means or of the time-averages, etc.). Even when the dynamics sequence converges itself, the equilibria to which it converges to (Nash equilibrium or social optimum) might not be unique \cite{7479108}. Approaching this question, Krichene et al. estimates the learning dynamics within an online learning model of player dynamics \cite{krichene2018tcps} and studies the traffic dynamics under partial control \cite{7604118}. 

To better understand the routing game in these systems and provide bigger impacts in transportation, we might also want to leverage concepts from other areas of engineering. Modeling and analysis of conflict, the study of systems, and the study of their equilibria are not unique to the routing game. In fact, very similar concepts and analysis shows up in control theory, which oftentimes, similarly to game theory, models interactions between a system of controllers (the "decision-makers") and their environment \cite{recht2019tour}. Unsurprisingly, many influences within game theory show up in many areas of controls (differential games, teams games, and distributed control) \cite{marden2018game}, and many connections can be made between the two. Even though there exists many similarities and influences between game theory and controls, there are many useful control design techniques that exists separate from game theory \cite{goodwin2001control} \cite{alessio2009survey} \cite{bansal2017hamilton}, but are conceptually similar/applicable to designing decision-makers within routing games, especially for routing systems like google. 

Therefore, in this article, we contribute to the literature of routing games and draw from other areas of engineering to make deeper connections and make use of control techniques for routing. We specifically introduce a new game design framework for routing games using the linear-quadratic regulator (LQR).




\subsection{Key Contributions}
The present article focuses on transforming the repeated game problem into a control theoretic problem, and studying the  convergence of the game dynamics in this framework. The key contributions of the article include the following:
\begin{itemize}
    \item We provide a  game design scheme for repeated routing games through the Linear Quadratic Regulator (LQR).
    \item We design control dynamics with specific choices of system parameters to respect the conservation of flow and achieve specific convergence results.
    \item We provide a method for producing piecewise affine optimal routing strategies using explicit model predictive control (MPC) techniques.
    \item We derive new theoretical results for specific cases of the new design scheme. 
    \item We illustrate a geometric framework for visualizing optimal routing solutions for every feasible state in the routing game.
    \item We extend the comparison literature between game theory and controls.
\end{itemize}
\subsection{Organization of the article} 
This article, first, explains the new control theoretic game design scheme in \textit{\nameref{Sec:LQR_Scheme_Sec}}, under  which we also build up the intuition and understanding of the framework by illustrating an evolution of the game from a one-shot framing (see \textit{\nameref{Subsec:One-shot}}) to a repeated game framing (see \textit{\nameref{Subsec:Repeated}}), and finally to the new control framing (see \textit{\nameref{Subsec:LQR}}). We further present analogies between the new framework and a classical algorithmic game theory model of routing games within this section (see \textit{\nameref{Subsec:LQR}}); then, we finally talk about how to design the game (see \textit{\nameref{Subsubsec:GameDesign}}) and its interpretation (see \textit{\nameref{Subsubsec:Interp}}). Next, we run through the system analysis and its properties in \textit{\nameref{Sec:SysAnalysis}}, where we additionally describe a few requirements for the system. In \textit{\nameref{Sec:MPC}}, we describe a algorithmic approach to solve the control routing problem using explicit MPC techniques. Finally, in \textit{\nameref{Sec:Results}} we run through a few numerical simulations and present its results with associated discussions.

\newpage

\section{Preliminary}

\begin{table}[h]
\footnotesize
\begin{tabular}{|p{1.7cm}||p{6.7cm}|p{6.3cm}|}
\hline
  \textbf{Notation}
  & \textbf{Description} 
  & \textbf{Type}\\
  \hline 
  
  $\network = (\vertexset, \edgeset)$ &A parallel network for the routing game made of a set of edges $\edgeset$ (as roads) and a set of vertices $\vertexset$ (one origin and one destination)&A directed graph\\
  \hline
  $\edgeelem \in \edgeset$&An edge in the edge set&An element of $\edgeset$\\
  \hline
  $\latfunc(\cdot)$&Affine latency function on edge $\edgeelem$
  &$\edgeflowset \mapsto \costset: \latfunc(\edgeflow) = a_e \edgeflow +b_e, a_e, b_e \geq 0$\\
  \hline 
  
  $\numedges$&   Number of parallel edges 
  &$\mathbb{N}$ \\
  \hline 

  $T$&   Time horizon (total number of rounds in game)
  & $\mathbb{N}$ (discrete time)\\
  \hline 
  $\flowset$&Probability Simplex
  & $\flowset:=\{ \flow: \sum_{\edgeelem \in \edgeset} \edgeflow = 1, \edgeflow \geq 0 \}$ \\
  \hline 
 $\flow$
  &Flow allocation vector on the network
  &  $\flowset, \flow = \left(\edgeflow\right)_{\edgeelem \in \edgeset}$\\
  \hline 
  $\loss$
  & Latency vector incurred on the network
  & $\lossvectorset,\loss = \left(\latfunc(\edgeflow)\right)_{\edgeelem \in \edgeset}$ \\
  \hline 
  $\stateiter$&System state vector at time $t$
  & $\flowset, \stateiter = (x_{1}^{\timesym},\dots, x_{\numedges}^{\timesym})^\top$\\
  \hline 
  $\controliter$
  &System input vector at time $t$
  & $\flowset,
\controliter = (u_{1}^\timesym,\dots, u_{\numedges}^{\timesym})^\top$\\
  \hline 
  $\state\optimal$&Steady state of the system&$\flowset$\\
  \hline
  $\initstate$&Initial state of the system&$\flowset$\\
  \hline
  $Q_f$
  &Terminal cost matrix
  & $\mathbb{R}^{\numedges\times \numedges}$ positive semi-definite\\
  \hline 
  $R$
  &Running cost matrix for inputs
  & $\mathbb{R}^{\numedges\times \numedges}$ positive semi-definite\\
  \hline 
  $Q$
  &Running cost matrix for states
  & $\mathbb{R}^{\numedges\times \numedges}$ positive semi-definite\\
  \hline 
  $\gamma$
  & Weight parameter in the dynamics
  & $\gamma \in [0,1]$\\
  \hline 
  $A$
  & Mapping from the memory to decision
  & $\mathbb{R}^{\numedges\times \numedges}$ left stochastic matrix\\
  \hline 
  $B$
  & Mapping from suggested routing to decision
  & $\mathbb{R}^{\numedges\times \numedges}$ left stochastic matrix\\
  \hline
\end{tabular}
\end{table}

\section{LQR Repeated Game Framework for Routing in Parallel Networks}
\label{Sec:LQR_Scheme_Sec}
We consider a specific routing game design problem where the routing game is posed with a parallel network, $\network$. $\network$ is modeled as a directed graph $\network = (\vertexset, \edgeset)$, where $\vertexset$ is a set containing two nodes (an origin $\origin$ and a destination $\destination$), and $\edgeset$ is an edge set containing $\numedges$ edges between each vertex.

\subsection{One-shot Game Framework}
\label{Subsec:One-shot}
We start by introducing the flow allocation problem, which we will use later as a potential equilibrium point for our repeated game. We will also refer to it as a one-shot game (i.e. a repeated game converging in one iteration).

\begin{definition}[Edge flows]
For each edge $\edgeelem \in \edgeset$, the edge flow, $\edgeflow \in \edgeflowset$, represents a large collection of non-atomic drivers using edge $\edgeelem$ on the network. Each edge flow, $\edgeflow$, represents a portion of traffic in the game (and thus has "mass"), while individual drivers themselves do not. They are negligible with respect to the total number of drivers present (i.e. the non-atomic framework). 
\label{Edge flows} 
\end{definition}

In the flow allocation problem, all drivers are allocated along the network's edges, and each allocation incurs a latency (i.e., loss), $\loss$ (this latency is usually labeled as the average travel time each driver took). Each edge's latency $\loss$ is modeled in the present work as an affine function associated with each edge.

\begin{definition}[Latency function] \label{Latency function} Each edge $\edgeelem \in \edgeset$ has an affine latency function $\latfunc(\cdot): \edgeflowset \mapsto \costset, \latfunc(\edgeflow) = a_e \edgeflow +b_e, a_e, b_e \geq 0$ that converts the edge flow $\edgeflow$ into the edge latency. The latency functions are nonnegative, nondecreasing and Lipschitz-continuous functions.
\label{def:latencyrepeated}
\end{definition}

Normally, in the allocation problem, drivers at the origin, $\origin$, are allocated strategically to reduce the total travel time the system will incur. The resulting allocation for all of the drivers of total mass $1$ in $\network$, is represented as a stochastic vector $\flow =  \left(\edgeflow\right)_{\edgeelem \in \edgeset} \in \flowvectorset$, where each $\edgeflow$ tells us the proportion of flow being allocated on edge $\edgeelem$. Specifically, 
\begin{linenomath*}
\begin{align}
    \flow \in \flowset\text{, where } \flowset := \{ \flow: \sum_{\edgeelem \in \edgeset} \edgeflow = 1, \edgeflow \geq 0 \} 
\end{align}
\end{linenomath*}

\subsubsection{Nash Equilibrium}
If we assume that all drivers act rationally and has access to perfect information, we can describe an equilibrium of the game, known as the \textit{\textbf{Nash equilibrium}}.

\begin{definition}[Nash equilibrium]\label{defi:nashrepeated}
A flow allocation $\flow \in \flowset$ is a Nash equilibrium (user equilibrium \cite{patriksson}, non-atomic equilibrium flow \cite{Nisam2007}, or Wardrop equilibrium \cite{Wardrop1952}) if and only if:
\begin{linenomath*}
\begin{equation}
    \forall \edgeelem \in \edgeset \text{ such that }\edgeflow > 0, ~~\latfunc(\edgeflow) =\min_{\otheredge \in \edgeset} \otheredgelatfunc(\otheredgeflow).
\end{equation}
\end{linenomath*}
\end{definition}

\subsection{Repeated Game Framework}
\label{Subsec:Repeated}
\label{section:RepeatedGames}
In order to encompass time-varying changes in the system, one can use "repeated play" or learning models \cite{cesa2006prediction}. In this framework, a flow allocator (i.e. Gooogle maps) makes a decision iteratively, as opposed to one time in the one-shot game, and uses the outcome of each iteration to adjust their next decision. Adding a $\timesym$ dependency for time to the aforementioned notation, $\flowiter \in \flowset$ becomes the flow allocation at time $\timesym$, $\lossiter \in \lossvectorset$ becomes a vector of latencies produced from functions $\latfunc(\cdot)$ at time $\timesym$, and throughout the game, the flow allocator iteratively chooses $\flowiter$ and then observes $\lossiter$. In this framework, designing learning models of player decisions and whether the resulting dynamics asymptotically converges to equilibrium is one of the objectives of our article. Furthermore, we are also interested in situations in which convergence happens in one shot. Previous research has been done to characterize classes of learning algorithms that converge to different sets of equilibria \cite{blum2006routing} \cite{marden2018game} \cite{hart2001general}. Specifically, different models of learning algorithms and their associated convergence guarantees has been studied for the routing game \cite{krichene2018tcps}, \cite{7604118}, \cite{krichene2015congestion}.



\subsection{LQR Framework}
\label{Subsec:LQR}
The following table sets the framework for the proposed work. The right column defines the repeated game framework previously described (using example designs from previous work \cite{NIPS2015_5843}), while the left column defines the classical control framework used later for LQR control. Both frameworks present obvious mathematical analogies \cite{marden2018game} \cite{marden2015game}.\\

\begin{table}[h]
\footnotesize

\begin{tabular}{|p{2.9cm}||p{6.1cm}|p{6.1cm}|}
 \hline
 &\textbf{Control problem} &  \textbf{Repeated game (online learning)}\\
 \hline
 Notation   & $\stateiter$ state, $\controliter$ control at time $\timesym$ & $\flowiter$ flow allocation, $\lossiter$ loss at time $\timesym$\\
  \hline
 Entities involved &\begin{itemize}
     \item Central decision maker: flow allocator
    \item Players: drivers
\end{itemize}& 
\begin{itemize}
    \item Central decision maker: flow allocator
    \item Players: flow allocator
\end{itemize}\\
\hline
 Design   & $\controliter(\cdot)$  s.t. $\stateiterplus = \ctrldyn(\stateiter,\controliter)$ is stable,  where $\ctrldyn(\cdot, \cdot)$ is the recurrence relation defining the dynamics.    &$\flowdyn(\cdot,\cdot)$ s.t. $\flowiterplus = \flowdyn(\flowiter,\lossiter)$ defines the dynamics for the routing game.\\
  \hline
  Example designs &$\controliter(\stateiter) = K\stateiter$ and $g(\stateiter,\controliter) = A\stateiter + B\controliter\implies \stateiterplus = (A+BK)\stateiter$, where $(A, B)$ is stabilizable, and $A+BK$ is asymptotically stable.& $\flowdyn(\flowiter,\lossiter) = \mathrm{argmin}_{\flowstoch \in \flowset} l_t^\top( \flowstoch - \flowiter)+\frac{1}{\eta_t}D_{\Psi}(\flowstoch, \flowiter)$, where $\eta_t$ is the learning rate at time \timesym, $D_{\Psi}(\cdot, \cdot)$ is the Bregman
divergence induced by a strongly
convex function $\Psi$ defined as $D_{\Psi}(x, y) = \Psi(x)-\Psi(y)-\nabla\Psi(y)^\top(x-y)$. This is called the \textbf{mirror descent algorithm} with Bregman divergence $D_{\Psi}$.\\
  \hline
 Objective function (to minimize)&Cumulative cost, $\sum_{t = 0}^{T-1} \costfunciter(\stateiter,\controliter) + \costfunciterT(\state_T) $&Cumulative regret, $\sum_{t = 0}^{T} \flowcostfunciter(\flowiter, \lossiter) - \min_{\flowstoch \in \flowset}\sum_{t = 0}^{T} \flowcostfunciter(\flowstoch, \lossiter)$\\
\hline
Example objective functions &$\costfunciter(\stateiter,\controliter) = \stateiter^\top Q \stateiter + \controliter^\top R \controliter$, 
$\costfunciterT(\stateiter) = \stateiter^\top Q_f \stateiter$, 
$Q,Q_f\succeq0,R\succ0$
& $\flowcostfunciter(\flowiter,\lossiter) = \flowiter^\top \lossiter$ \\
\hline
\end{tabular}
\end{table}

We re-frame the routing problem under the repeated game framework into the LQR framework and solve for/analyze player strategies using existing control techniques. We will essentially aim to do the same game theoretical analysis within linear/nonlinear controls. In the right column of the table, we list the setup under the repeated game framework. The goal of the routing problem there is to design algorithms that dictate strategies for players (the flow allocator). At each time step $\timesym$, the algorithm updates the flow allocation $\flowiterplus$ based on the players' memory $\flowiter$ and the loss on the edge $\lossiter$ induced by the traffic \cite{krichene2015congestion} \cite{7604118} \cite{pmlr-v37-krichene15} \cite{roughgarden2010algorithmic}. The goal is to guarantee a sub-linear regret and a convergence to the set
of Nash equilibria.\\

We start from the repeated games setup, and cast the repeated game problem in a control theoretic framework. We keep the authoritative central decision maker $\control$ as the flow allocator, but we differentiate between the "player" and the "decision-maker" as coined in the original game theoretic setup. Here, the "players" will refer to the drivers in the network flow being routed. In our framework, the "players" do not play; they are only considered by the routing system and the model. The central decision maker is what we want to design in this framework. It has full access to the flow allocation $\state$ on the network, and aims to achieve a certain routing goal (e.g., steer the network flow to a target flow distribution, or minimize the average travel time). The central decision maker embeds the routing goal in the cumulative cost function $\sum_{t = 0}^{T-1} \costfunciter(\stateiter,\controliter) + \costfunciterT(\state_T)$. The players now update their routing decision $\stateiterplus$ by leveraging their memory $\stateiter$ and the central decision maker's suggested routing distribution $\controliter$ at each time step. The target flow allocations $\state\optimal$ are characterized as the steady states of the controlled system, and the central decision maker chooses the optimal controls $\{\controliter\}$  to minimize the cumulative cost function and stabilizes the system to its steady states.\\

The choice of LQR is motivated by the observation that a potential game with affine latency and a parallel network can be naturally framed as a \textbf{quadratic} cost and \textbf{linear} dynamic problem, leveraging the convex formulation of the Rosenthal potential:
\begin{linenomath*}
\begin{align}
    \label{sys:lqr}
    &\underset{\controliter, t = 0, \dots, T-1}{\text{minimize}}  &\sum_{t = 0}^{T-1} \left(\stateiter^\top Q \stateiter + \controliter^\top R \controliter\right)
    &+ \state_T^\top Q_f\state_T \\
    &\text{subject to}&
    \stateiterplus &= A\stateiter + B\controliter, t = 0, 1, \dots, T-1\\
    &&\initstate &= \state(0)
\end{align}
\end{linenomath*}
    where $\stateiter \in \statevectorset$ is the state of the system at time $\timesym$ and $\controliter \in \controlvectorset$ is the action of the decision-maker at time $\timesym$, with appropriate matrices $A, B, Q, R, Q_f$. 
Within this framework, we can leverage known optimal control techniques to formulate optimal strategies and inherent known properties for the system. Furthermore, definition \eqref{defi:nashrepeated} for the context here becomes:

\subsubsection{Nash Equilibrium}
\begin{definition}[Nash equilibrium]\label{defi:nash_ctrls}
A state of the system $\state \in \flowset$ is a Nash equilibrium if and only if:
\begin{linenomath*}
\begin{equation}
    \forall \edgeelem \in \edgeset \text{ such that }\edgestate > 0, ~~\latfunc(\edgestate) =\min_{\otheredge \in \edgeset} \otheredgestatelatfunc(\otheredgestate).
\end{equation}
\end{linenomath*}
\end{definition}

We reiterate the key differences between the repeated game framework and our control theoretic framework:

\begin{enumerate}
    \item In repeated games for parallel networks, the two entities involved, the player and the central decision maker are the same (a flow allocator). However, in our control theoretic framework, the two entities are different -- \textbf{players being the drivers} (in the present case aggregated as non-atomic along the edges) and \textbf{the central decision-maker being the flow allocator}.
    \item In repeated games, the players update their flow allocation based on their memory $\flowiter$ and the incurred loss $\lossiter$ on the network. In our control theoretic framework, the players update their flow allocation based on \textbf{their memory} $\stateiter$ and \textbf{the central decision maker's suggested routing decision $\controliter$}.
    \item In repeated games, one possible goal of the non-atomic player, the flow allocator, (see for example \cite{pmlr-v37-krichene15} \cite{krichene2015heterogeneous}) is to update their flow allocation in order to converge to the set of Nash equilibria and to achieve sub-linear regret. In our control theoretic framework, \textbf{the goal of the central decision maker is to steer the flow allocation to a target flow distribution, by inputting suggestions to the system's players (drivers)}.
\end{enumerate}

\subsubsection{Designing Games with LQR}
\label{Subsubsec:GameDesign}
We now describe our control routing model.\\
\begin{enumerate}
    \item States and controls: the \textit{state} of the LQR at time $\timesym$ is the flow vector \stateiter, consisting of the flow on each edge at time $t$, and the \textit{control} at time $\timesym$ is the flow allocation $ \controliter$ given by the decision maker, i.e.,
    \begin{linenomath*}
    \begin{align}
    \stateiter = (x_{1}^{\timesym},\dots, x_{\numedges}^{\timesym})^\top&, t = 0, \dots, T\\
    \controliter = (u_{1}^{\timesym},\dots, u_{\numedges}^{\timesym})^\top&, t = 0, \dots, T-1.
    \end{align}
    \end{linenomath*}
    
    We model both $\stateiter,\controliter \in \flowset$ at each time step \timesym.
    \item Equation for the dynamics: we assume the dynamics of the flow update is described using a linear time invariant difference equation
    \begin{linenomath*}
    \begin{align}
    \stateiterplus = \gamma A \stateiter + (1-\gamma)B \controliter, t = 0, 1, \dots, T-1
    \end{align}
    \end{linenomath*}
    where $\gamma \in [0,1]$ is a design parameter that weights the contribution of $\stateiter$ and $\controliter$ in the update rule, and $A, B \in \mathbb{R}^{\numedges \times \numedges}$ are two design matrices deciding how the flow vectors $\stateiter$ and $\controliter$ will affect the drivers' routing decision in the subsequent time step respectively.
    \item Cost function: the cost function of our control problem is a summation of quadratic functions
    \begin{linenomath*}
    \begin{align}
        \sum_{t = 0}^{T-1} \left(\stateiter^\top Q \stateiter + \controliter^\top R \controliter \right)+ \state_T^\top Q_f \state_T
    \end{align}
    \end{linenomath*}
    where $Q, Q_f, R \succeq 0$.
\end{enumerate}


The resulting control problem is
\begin{linenomath*}
\begin{align}  
    &\underset{\controliter, t = 0, \dots, T-1}{\text{minimize}} &\sum_{t = 0}^{T-1} \left(\stateiter^\top Q \stateiter + \controliter^\top R \controliter \right)+ \state_T^\top Q_f \state_T& \label{begin_lqrrouting}\\ 
    &\text{subject to }&\stateiterplus = \gamma A \stateiter + (1-\gamma)B \controliter&, t = 0, 1, \dots, T-1\\
    & & \stateiter \in \flowset&, t = 0, \dots, T \\
    & & \controliter \in \flowset&, t = 0, \dots, T-1 \label{end_lqrrouting}
\end{align}
\end{linenomath*}

\subsubsection{Interpretation of the Framework}
\label{Subsubsec:Interp}
We now provide interpretations for the design of the control system. 
\begin{enumerate}
    \item States $\stateiter$: 
    the interpretation of the state $\mathbf{x}_t$ is twofold: (1) one could interpret it as the routing decision made by the drivers from one day to another, and (2) it is also the \textit{actual flow distribution} on the network, incurred by the drivers. We assume all drivers comply with a routing decision which depends at the same time on the past state, and the control.
    \item Controls $\controliter$: One could interpret the control $\controliter$ as the recommended routing decision from the central decision maker.
    \item The $A\stateiter$ term in the dynamics represents the mapping from the drivers' memory (past state inherited from the day before) to their  flow allocation decision. 
    \item The $B\controliter$ term in the dynamics represents the mapping from the central decision maker's suggested routing distribution to the players' flow allocation decision.
    \item When updating their flow allocation decision, the drivers leverage their memory $\stateiter$ and the suggested routing distribution $\controliter$ from the decision maker. $\gamma\in [0,1]$ is the design parameter that characterizes the trade-off between the two contributions.
    \item (Quadratic) Cost function: the design of the quadratic cost function $\stateiter^\top Q \stateiter + \controliter^\top R \controliter$ is motivated by the observation that the Rosenthal potential function and the cost function associated to social welfare in a potential game are both quadratic for linear latency functions \cite{Rosenthal1973ACO}.\\
    
    \begin{definition}[Rosenthal potential] The Rosenthal potential function of a potential game, with a flow vector $\state$, is defined as
    \begin{linenomath*}
    \begin{align}
        \mathcal{J}_1(\state) = \sum_{\edgeelem \in \edgeset}\int_0^{x_e} \latfunc(y)\,\mathrm{d}y.
    \end{align}
    \end{linenomath*}
    \end{definition}
    
    \begin{definition}[Social Welfare] The cost function associated to social welfare is defined as
    \begin{linenomath*}
    \begin{align*}
        \mathcal{J}_2(\state) = \sum_{\edgeelem\in \edgeset}x_e \cdot \latfunc(x_e),
    \end{align*}
    \end{linenomath*}
    \end{definition}
    
    The importance of the Rosenthal potential and the social welfare is that the Nash equilibrium and social optimal of the routing game can be characterized as the local minimizers of these functions respectively, respecting the conservation of flow. With linear latency parallel networks,
    \begin{linenomath*}
    \begin{align}
        \mathcal{J}_1(\state) = \sum_{i=1}^{\numedges}\int_0^{x_i} (a_i y+b_i)\,\mathrm{d}y = \sum_{i=1}^{\numedges} \frac{a_ix_i^2}{2} +b_i x_i = \state^\top \tilde{Q} \state + \mathbf{b}^\top \state,
    \end{align}
    \end{linenomath*}
where $\tilde{Q} = \text{diag}(\frac{a_1}{2},\dots, \frac{a_{\numedges}}{2}), \mathbf{b} = (b_1, \dots, b_{\numedges})$ and 
\begin{linenomath*}
    \begin{align}
        \mathcal{J}_2(\mathbf{x}) = \sum_{i = 1}^{\numedges} x_i (a_i x_i +b_i) = \sum_{i=1}^{\numedges} a_ix_i^2 +b_i x_i = 2\state^\top \tilde{Q} \state + \mathbf{b}^\top \state,
    \end{align}
    \end{linenomath*}
    with the $\tilde{Q}, \mathbf{b}$ defined as the same as above. Note that $\tilde{Q}$ and \textbf{b} are used to construct $Q$ and $Q_f$ of the cost function.

    \item Design of the matrix $A$ in the dynamics: since the $A\mathbf{x}_t$ is the mapping from the memory to the flow allocation decision, the matrix $A$ should reflect how the drivers take the memory into account when updating their routing decisions.
    \item Design of the matrix $B$ in the dynamics: since the $B\mathbf{u}_t$ is the mapping from the suggested routing distribution to the flow allocation decision, the matrix $B$ should reflect how the drivers take the central decision maker's suggested flow into account when updating their routing decisions.
\end{enumerate}

\section{System Analysis and System Properties} \label{Sec:SysAnalysis}
One of the contributions of the work is about how to design specific A and B matrices to produce various convergence results.
\begin{enumerate}
    \item Construction of $A, B$: when designing the matrices $A, B$ in the dynamics, we need to take the conservation of flow into account. More specifically, $A, B$ have to be chosen such that the following property holds
    \begin{linenomath*}
    \begin{equation}
        \stateiter, \controliter \in \flowset \quad \implies  \quad \gamma A\stateiter + (1-\gamma)B\controliter = \stateiterplus \in \flowset.
        \label{constraints:flow_cons}
    \end{equation}
    \end{linenomath*}
    One way to ensure the conservation of flow on the network is to construct \textit{left stochastic} $A, B$.\\
    \begin{definition}
    A matrix $P=(\mathbf{p}_1, \dots, \mathbf{p}_{\numedges})\in \mathbb{R}^{\numedges\times \numedges}$ is \textit{left stochastic} if $\mathbf{p}_i \in \Delta, i = 1, \dots, \numedges$.
    \end{definition}
    \begin{lemma}If $A, B\in \mathbb{R}^{\numedges\times \numedges}$ are left stochastic, and $\initstate, \controliter \in\flowset, t = 0, \dots, T-1 $, then $\gamma A\stateiter +(1-\gamma) B\controliter\ = \stateiterplus \in \flowset, t = 0, \dots, T-1$.
    \begin{proof}
    We will show that if $\stateiter, \controliter \in \flowset$, then $\gamma A\stateiter +(1-\gamma) B\controliter = \stateiterplus \in \flowset$. The proof is then complete by induction on $t$. \\
    First observe that $\mathbf{x}_{t+1}$ has non-negative entries due to the non-negativity of the entries in $A, B, \stateiter, \controliter$. 
    Furthermore,
    \begin{linenomath*}
    \begin{align}
        \mathbb{I}^\top \stateiterplus &= \gamma \mathbb{I}^\top A \stateiter + (1-\gamma)\mathbb{I}^\top B \controliter\\
        &= \gamma \mathbb{I}^\top \stateiter + (1-\gamma)\mathbb{I}^\top \controliter\\
        & = \gamma + (1-\gamma) = 1.
    \end{align}
    \end{linenomath*}
    \end{proof}
    \label{lemma:flow_cons}
    \end{lemma}
    \item Examples: the analysis of the LQR problem with probability simplex constraints is in general not tractable. In other words, we are not aware of a way to solve the discrete-time  algebraic  Ricatti equation in a way that preserves \eqref{constraints:flow_cons}. We therefore provide the important special cases where either (1) exact calculations of the steady state and optimal controls are allowed, or (2) interesting behavioral interpretation of the players can be given.
    \begin{example}
    We first consider the special case where
    \begin{equation}
    B=\begin{pmatrix}
        \frac{1}{\numedges}& \frac{1}{\numedges}&\dots &\frac{1}{\numedges}\\
        \frac{1}{\numedges}& \frac{1}{\numedges}&\dots &\frac{1}{\numedges}\\
        \vdots  & \vdots  & \ddots & \vdots  \\
        \frac{1}{\numedges}& \frac{1}{\numedges}&\dots &\frac{1}{\numedges}\\
        \end{pmatrix}\in \mathbb{R}^{\numedges\times \numedges}.
    \end{equation}
        This choice of $B$ matrix in the dynamics can be interpreted as follows: when the players (drivers) leverage the central decision maker's suggested routing distribution $\controliter$ to update their flow allocation on each edge, they simply average the suggested flow at each time step $\timesym$, and take only the average $\frac{1}{\numedges}\sum_{i = 1}^{\numedges} u_i^t$ into account.\\
        We show that when the players (drivers) update their flow allocation in this way, the controller $\controliter$ is essentially ``muted'', and therefore a target state may not be controllable. The state $\stateiter$ will converge to the same equilibrium that is independent of the choice of the cost function.
        \begin{property}
            Consider the control problem \eqref{begin_lqrrouting}-\eqref{end_lqrrouting} with 
            \begin{equation}\label{eq:1_3_B}
               B =
               \begin{pmatrix}
        \frac{1}{\numedges}& \frac{1}{\numedges}&\dots &\frac{1}{\numedges}\\
        \frac{1}{\numedges}& \frac{1}{\numedges}&\dots &\frac{1}{\numedges}\\
        \vdots  & \vdots  & \ddots & \vdots  \\
        \frac{1}{\numedges}& \frac{1}{\numedges}&\dots &\frac{1}{\numedges}\\
        \end{pmatrix},
        \gamma \in [0,1), \text{ and any left stochastic matrix }A, 
        \end{equation}
        then the unique steady state of the system is given by
        \begin{equation}
            \state\optimal = \frac{1-\gamma}{\numedges}(I_{\numedges}-\gamma A)^{-1}\mathbb{I},
        \end{equation}
        where $I_{\numedges}$ is the $\numedges$-dimensional identity matrix, and $\mathbb{I} = (1, 1, \dots, 1)^\top \in \mathbb{R}^{\numedges}$ is the all-ones vector.
        \begin{proof}
        The key observation is, with equation (\ref{eq:1_3_B}), since we have $\controliter \in \flowset, t = 0, \dots, T-1$, 
        \begin{linenomath*}
        \begin{align}
            B\controliter = \frac{1}{\numedges} \mathbb{I},
        \end{align}
        \end{linenomath*}
        therefore
        \begin{linenomath*}
        \begin{align}
            \stateiterplus=\gamma A\stateiter + (1-\gamma)\frac{1}{\numedges} \mathbb{I},
        \end{align}
        \end{linenomath*}
        and the steady state should satisfy
        \begin{equation}
            \state\optimal=\gamma A\state\optimal + (1-\gamma)\frac{1}{\numedges} \mathbb{I}, \quad \implies \quad \state\optimal = \frac{1-\gamma}{\numedges}(I_{\numedges}-\gamma A)^{-1}\mathbb{I}, 
        \end{equation}
        with the invertibility of $I_{\numedges}-\gamma A$ guaranteed by the fact that the spectral radius of any left stochastic matrix is at most $1$, and $\gamma<1$.
        \end{proof}
        \end{property}
        The foregoing property suggests that in the case of equation \eqref{eq:1_3_B}, the flow allocation on the network does converge to an equilibrium, but the equilibrium depends only on the design parameters $A, \gamma$ and the number of edges $\numedges$, which can be different from the Nash equilibrium.
        \end{example}
        \begin{example}
        We next consider the special case where
        \begin{equation}
        A = \begin{pmatrix}
        \frac{1}{\numedges}& \frac{1}{\numedges}&\dots &\frac{1}{\numedges}\\
        \frac{1}{\numedges}& \frac{1}{\numedges}&\dots &\frac{1}{\numedges}\\
        \vdots  & \vdots  & \ddots & \vdots  \\
        \frac{1}{\numedges}& \frac{1}{\numedges}&\dots &\frac{1}{\numedges}\\
        \end{pmatrix}\in \mathbb{R}^{\numedges\times \numedges}.
        \end{equation}
        This choice of $A$ matrix in the dynamics can be interpreted as follows: when the players (drivers) leverage 
        their memory $\stateiter$ to update their flow allocation on each edge, they simply average the past flow at each time step $\timesym$, and take only this average $\frac{1}{\numedges}\sum_{i = 1}^{\numedges} x_i^t$ into account.\\
        We similarly show that when the players (drivers) update their flow allocation in this way, their memory has no impact in the update rule.
        \begin{property} \label{1/N_A}
            Consider the control problem \eqref{begin_lqrrouting}-\eqref{end_lqrrouting} with 
            \begin{equation}\label{eq:1_3_A}
               A =
               \begin{pmatrix}
        \frac{1}{\numedges}& \frac{1}{\numedges}&\dots &\frac{1}{\numedges}\\
        \frac{1}{\numedges}& \frac{1}{\numedges}&\dots &\frac{1}{\numedges}\\
        \vdots  & \vdots  & \ddots & \vdots  \\
        \frac{1}{\numedges}& \frac{1}{\numedges}&\dots &\frac{1}{\numedges}\\
        \end{pmatrix} \text{ and any left stochastic matrix } B,  
        \end{equation}
        then the optimal controller $\controliter, t = 0, \dots, T-2$ is the solution set of the optimization problem
        \begin{linenomath*}
        \begin{align}
    \underset{\control}{\text{minimize}}\quad&\left( \frac{\gamma}{\numedges} \mathbb{I} + (1-\gamma)B\control\right)^\top Q \left( \frac{\gamma}{\numedges} \mathbb{I} + (1-\gamma)B\control\right) + \control^\top R \control\\
    \text{subject to}\quad& \control \in \flowset,
\end{align}
\end{linenomath*}
and the optimal $\control_{T-1}$ is the solution set of the optimization problem
\begin{linenomath*}
\begin{align}
    \underset{\control}{\text{minimize}}\quad&\left( \frac{\gamma}{\numedges} \mathbb{I} + (1-\gamma)B\control\right)^\top Q_f \left( \frac{\gamma}{\numedges} \mathbb{I} + (1-\gamma)B\control\right) + \control^\top R \control\\
    \text{subject to}\quad& \control\in \flowset.
\end{align}
\end{linenomath*}
        \begin{proof}
        The key observation is, with equation (\ref{eq:1_3_A}), since we have $\stateiter \in \flowset, t = 0, \dots, T-1$, 
        \begin{linenomath*}
        \begin{align}
            A\stateiter = \frac{1}{\numedges} \mathbb{I},
        \end{align}
        \end{linenomath*}
        therefore
        \begin{linenomath*}
        \begin{align}
            \stateiterplus=\gamma \frac{1}{\numedges} \mathbb{I} + (1-\gamma)B\controliter.
        \end{align}
        \end{linenomath*}
        Therefore, the optimal control of problem  \eqref{begin_lqrrouting}-\eqref{end_lqrrouting} can be explicitly solved by the optimization problem
        \begin{linenomath*}
        \begin{align}
    \underset{\controliter, t=0,\dots T-1}{\text{minimize}}\quad&\initstate^\top Q \initstate+ \sum_{t = 0}^{T-2} \left(\left( \frac{\gamma}{\numedges} \mathbb{I} + (1-\gamma)B\controliter\right)^\top Q \left( \frac{\gamma}{\numedges} \mathbb{I} + (1-\gamma)B\controliter\right) + \controliter^\top R \controliter \right)+\nonumber\\
    &\left( \left(\frac{\gamma}{\numedges} \mathbb{I} + (1-\gamma)B\mathbf{u}_{T-1}\right)^\top Q_f \left( \frac{\gamma}{\numedges} \mathbb{I} + (1-\gamma)B\control_{T-1}
    \right) + \control_{T-1}^\top R \control_{T-1}\right)\\
    \text{subject to}\quad& \controliter \in \flowset, t = 0, \dots, T-1, 
\end{align}
\end{linenomath*}
which is further equivalent to solving the constrained optimization problems
\begin{linenomath*}
\begin{align}
    \underset{\control}{\text{minimize}}\quad&\left( \frac{\gamma}{\numedges} \mathbb{I} + (1-\gamma)B\control\right)^\top Q \left( \frac{\gamma}{\numedges} \mathbb{I} + (1-\gamma)B\control\right) + \control^\top R \control\\
    \text{subject to}\quad& \control \in \flowset
\end{align}
\end{linenomath*}
and
\begin{linenomath*}
\begin{align}
    \underset{\control}{\text{minimize}}\quad&\left( \frac{\gamma}{\numedges} \mathbb{I} + (1-\gamma)B\control\right)^\top Q_f \left( \frac{\gamma}{\numedges} \mathbb{I} + (1-\gamma)B\control\right) + \control^\top R \control\\
    \text{subject to}\quad& \control \in \flowset.
\end{align}
\end{linenomath*}
If $Q_f=Q$ and the solution of the above optimization problem exists, the minimizer $\control\optimal$ characterizes the steady state $\state\optimal = \frac{\gamma}{\numedges} \mathbb{I} + (1-\gamma)B\control\optimal$.
        \end{proof}
        \end{property}
        Note that in the very specific case where $R = 0$ and $Q, Q_f$ are defined by the Rosenthal potential function, the minimizer of the above optimization problem and steady state of the system $\state\optimal$ is the Nash equilibrium. If the linear equation system $\state\optimal = \frac{\gamma}{\numedges} \mathbb{I} + (1-\gamma)B\control$ has at least one solution $\control\optimal$ (e.g., when $B$ is invertible and $\gamma\neq 1$), then the  optimal controllers are $\controliter = \control\optimal, t = 0, \dots, T-1$. 
         \end{example}
         \begin{example}
         The third special case we would like to analyze is $A = B = I_{\numedges}$. The players update their flow allocation with $\stateiterplus = \gamma \stateiter + (1-\gamma) \controliter$ at each time step $\timesym$, which is the same as the element-wise update $x_i^{t+1} = \gamma x_i^t + (1-\gamma) u_i^t, i = 1, \dots, N, t = 0, \dots, T-1$.\\
         The interpretation of this choice of $A,B$ is that the players update the flow on each edge \textit{independently}: to update the flow on edge $i$, the players leverage their memory and the suggested routing decision on the edge $i$ only.
         \end{example}
\end{enumerate}


\section{Solutions via Explicit MPC: A Multiparametric Quadratic Programming Approach}
\label{Sec:MPC}
\begin{figure}[h]
\centering
\includegraphics[scale=0.2]{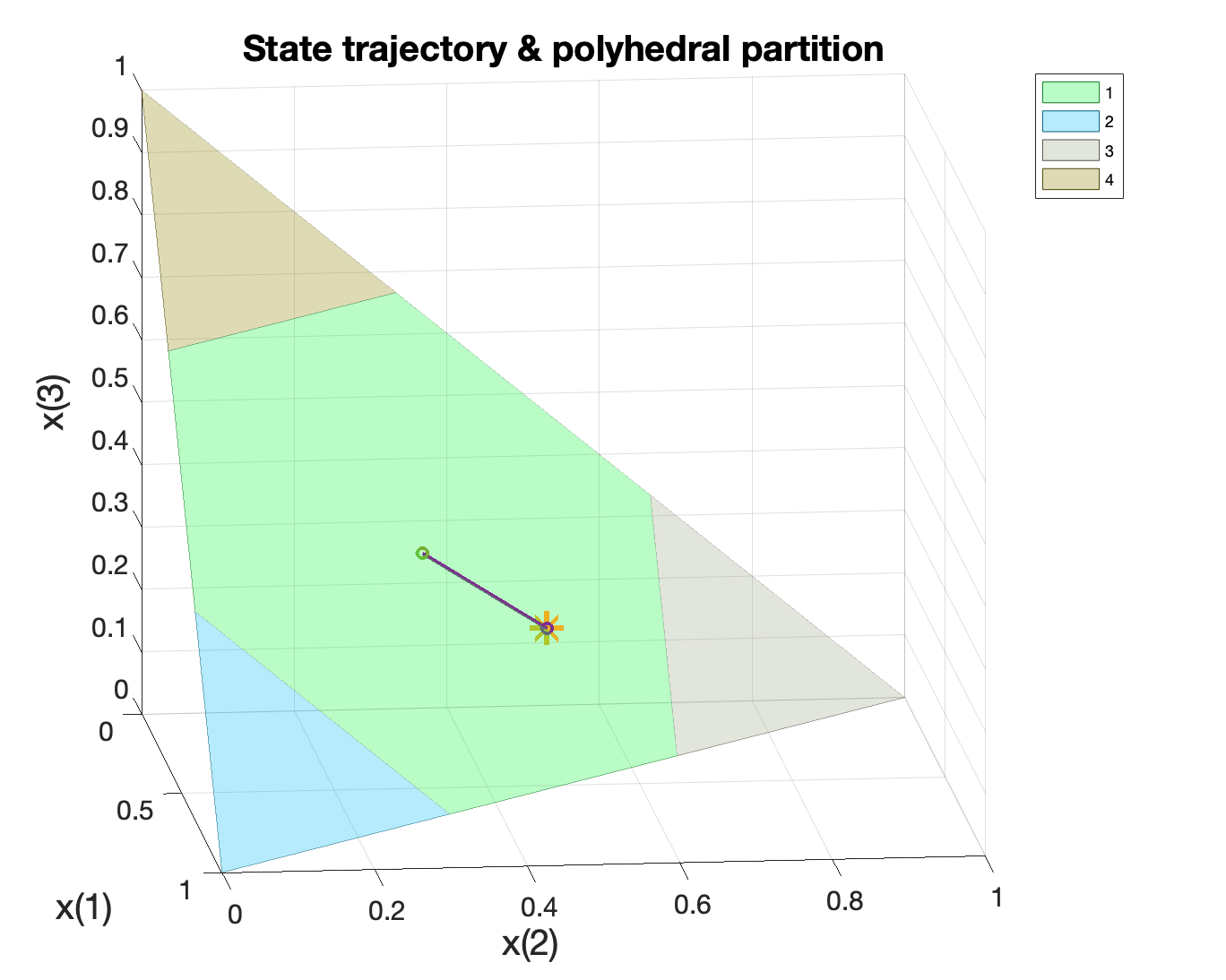}
\caption{Here, the optimal strategy for our player is space-varying and we visualize solutions with respect to the feasible state space of the routing game. Each colored region of the graph represents a unique affine strategy that is optimal for those points in space. For example, region 1 (in green) has the solution $\controliter\optimal(\stateiter)$ = $\begin{pmatrix} -\frac{2}{3} & \frac{1}{3} & \frac{1}{3}\\ \frac{1}{3} & -\frac{2}{3} & \frac{1}{3} \\
\frac{1}{3} & \frac{1}{3} & -\frac{2}{3} \\
\end{pmatrix} \stateiter + \begin{pmatrix} \frac{1}{3} \\
\frac{1}{3} \\
\frac{1}{3} \end{pmatrix}$.}
\label{fig:example_cr}
\end{figure}

At first glance, one might consider using the discrete-time algebraic Ricatti equation to obtain solutions for system \eqref{begin_lqrrouting}-\eqref{end_lqrrouting}. 
\begin{definition}[Discrete-time Algebraic Ricatti Equation] The discrete-time algebraic Ricatti equation is a nonlinear equation that gives a solution to the LQR problem presented in equation \eqref{sys:lqr}, described as $P_t$ evolving backwards in time from $P_T = Q_f$ according to
\begin{equation}
    P_{t-1} = Q + A\transpose P_{t} A - A\transpose P_{t} B(R + B\transpose Q_f B)\inv B\transpose P_{t} A.
    \end{equation}
\end{definition}

However, using the algebraic Riccati equation is not possible in the presence of our flow conservation constraints, as specified in equation \eqref{constraints:flow_cons}. According to lemma \eqref{lemma:flow_cons}, in order to design a flow preserving scheme and guarantee flow conservation, we must be able to guarantee that each $\stateiter$ and $\controliter$ vector given by a chosen controller remains on the probability simplex $\flowset$, while the optimal controller given by the algebraic Riccati equation drives the system to the steady state $\state\optimal = 0$. Therefore, we turn to multiparametric quadratic programming (mpQP), a technique used in explicit MPC to generate explicit piecewise affine controllers. To generate these explicit strategies, we will formulate the LQR problem into a mpQP problem, utilize an mpQP solver proposed in  \cite{bemporad2015multiparametric} to generate piecewise affine functions for $\controliter$ (the player needs to execute different affine controllers depending on the state), and use (within the solver) a geometric algorithm to plot regions, known as \textit{\textbf{critical regions}}, on the state space to visualize when players should use which specific optimal strategies. Optimal strategies produced by this technique will be piecewisely-affine with respect to the state of the game, hence, optimal strategies in this context are space-varying (see figure \eqref{fig:example_cr} for an example).

Algorithm \eqref{alg:mpc_algo} defines the process by which we use mpQP to generate solutions for routing games. Specifically, we consider the optimal control formulation in problem \eqref{begin_lqrrouting}-\eqref{end_lqrrouting}. We take the linear flow and input constraints in \eqref{constraints:flow_cons} and construct them in the form
\begin{linenomath*}
    \begin{align}
       A_u \controliter \leq b_u, \quad A_x \stateiter \leq b_x.
    \end{align}
    \end{linenomath*}
    To solve for the optimal controls $\control_0,\dots, \control_{T-1}$, first we condense the problem into a strictly convex quadratic programming problem
    \begin{linenomath*}
\begin{align}
    \label{sys:begincondenseqp}
    \condenseqpobj\optimal(x) &\triangleq& &\textit{minimize}& &\frac{1}{2} z^\top Pz + (Fx + c)^\top z + \frac{1}{2} x^\top Yx& \\
    &&&\textit{subject to}& & Gz \leq W + Sx&
    \label{sys:endcondenseqp}
\end{align}
\end{linenomath*}
and then use the mpQP algorithm from reference \cite{bemporad2015multiparametric} to solve for $z$ where $z := \left( \control_0^\top \dots \control_{T-1}^\top \right)^\top \in \R^{\nmpc}$ is the decision  variable.

$P=P^\top \in \R^{\nmpc \by \nmpc}$, $F \in \R^{\nmpc \by \mmpc}$, $c \in \R^{\nmpc}$, $Y \in \R^{\mmpc \by \mmpc}$ are defined by the cost matrices, and $G \in \R^{\nmpc \by \qmpc}$, $ W \in \R^{\qmpc}$, $S \in \R^{\qmpc \by \mmpc}$ define the constraints imposed by \eqref{begin_lqrrouting}-\eqref{end_lqrrouting} in compact form. The following property guarantees the existence of a (piecewise affine) solution of the condensed problem \eqref{sys:begincondenseqp}-\eqref{sys:endcondenseqp}.

\begin{property}[Existence of a piecewise affine solution \cite{bemporad2015multiparametric}] Consider \eqref{sys:begincondenseqp}-\eqref{sys:endcondenseqp} with $P=P^\top \succ 0$ and let $x\in \stateset=\R^n$.
\begin{enumerate}
    \item The set $\stateset_f$ of states for which the problem is feasible is a polyhedron.
    \item The optimizer function $z\optimal: \stateset_f \to \mathbb{R}^n$ is piecewise affine and continuous over $\stateset_f$.
    \item If 
   $ \begin{pmatrix}
        P &\F^\top\\
        \F &\Y \\
    \end{pmatrix} \succeq 0$ and symmetric, then the value function $\condenseqpobj\optimal$ is continuous, convex, and piecewise quadratic over $\stateset_f$.
\end{enumerate}
\end{property}

\begin{definition}[Critical Region \cite{spjotvold2006facet}] A critical region is a polyhedron in $\stateset$ for which there exists an optimal solution $z\optimal(\cdot)$ that is affine and optimal for the entire region. Each critical region is defined by a unique set of active constraints $\activeset$ that is common for all points in the region. 
\end{definition}

\IncMargin{1em}
\begin{algorithm}
    
    \SetKwInOut{Input}{Input}\SetKwInOut{Output}{Output}
    \SetKwData{qp}{$P,F,c,Y,G,W,S$}\SetKwData{input}{$Q$,$R$,$Q_f$,$N_c$,$A$,$B$,$A_u$,$b_u$,$A_x$,$b_x$}
    \SetKwFunction{Condense}{Condense}\SetKwFunction{mpqpAlgo}{MPQPSolver}
    \Input{\input}
    \Output{Optimal strategy solutions: \mpcctrl, \mpcctrladd\\ Critical Regions to apply solutions: \CR, \CRCheckH, \CRCheckK, }
    \BlankLine
    \qp $\leftarrow$ \Condense{\input}\;
    $\textbf{return}$ \mpqpAlgo{\qp}\;
    \caption{disjoint decomposition}\label{algo_disjdecomp}
    \caption{mpQP algorithm}\label{alg:mpc_algo}
\end{algorithm}\DecMargin{1em}

As explained in (\cite{bemporad2015multiparametric}), $\mpqpAlgo(\cdot)$ uses a geometrical algorithm (see \cite{spjotvold2006facet}) to find the critical regions of the feasible state space by understanding when sets of active constraints differs from point to point in the state space. Once it has identified a critical region, the algorithm is able to uncover the optimal solution defined within the region by using KKT conditions (or similar methods) to obtain the piecewise affine solution $z\optimal$. Algorithm \eqref{alg:mpc_algo} produces (1) the critical regions $\CR$, (2) the matrices to check if a state is in specific critical region, $\CRCheckH = \{\CRCheckH_r\}$ and $\CRCheckK = \{\CRCheckK_r\}$, and (3) the matrices that characterizes the affine solution $\mpcctrl = \{\mpcctrl_r\}$ and  $\mpcctrladd = \{\mpcctrladd_r\}$. Therefore, during the routing game, to recover the optimal strategies, the central decision maker executes
\begin{linenomath*}
\begin{align}
    \text{if } \CRCheckH_{\regionelem}\stateiter \leq \CRCheckK_{\regionelem} ~ \text{, then } \controliter\optimal = \mpcctrl_{\regionelem}\stateiter + \mpcctrladd_{\regionelem}
\end{align}
\end{linenomath*}
and \stateiter~belongs to critical region \regionelem. 





\section{Illustration on Simple Numerics}
\label{Sec:Results}
In this section, we plan to provide numerical results and discuss
\begin{enumerate}
    \item differences in varying $A$ matrices in the dynamics, and the interpretations,
    \item changes in $\gamma$ and its effect on the converging rate, 
    \item different initial states and its effect on the converging rate,
    \item different cost functions to drive the system towards different stable equilibrium.
\end{enumerate}

More specifically, we will compare between
\begin{enumerate}
    \item $A = \begin{pmatrix}
        \frac{1}{3}& \frac{1}{3}&\frac{1}{3}\\
        \frac{1}{3}& \frac{1}{3}&\frac{1}{3}\\
        \frac{1}{3}& \frac{1}{3}&\frac{1}{3}
        \end{pmatrix}$ and  $A = I_3$,  with all the other parameters fixed as $B = I_3,  Q = Q_f = I_3, R = 0, \gamma = 0.5, \initstate = (0.3, 0.5, 0.2)^\top$. This means we investigate the differences between the players' routing decision process when they ``mute'' their memory, and when they leverage their memory on each edge independently (i.e., when $x_i^{t+1}$ depends only on $x_i^t$ in the past state $\stateiter$ for $i = 1, \dots, \numedges$). 
    \item $\gamma = 0.5$ and $\gamma = 0.7$:  with all the other parameters fixed as $A = B = I_3, Q = Q_f = \begin{pmatrix}
         1 & 0 & 0 \\
         0 & 2 & 0 \\
         0 & 0 & 4
    \end{pmatrix}, R = \begin{pmatrix}
         0 & 0 & 0 \\
         0 & 0 & 0 \\
         0 & 0 & 0
    \end{pmatrix}, \mathbf{x}_0 = (0.3, 0.5, 0.2)^\top$.  This means we investigate the differences between the players' routing decision process when they weight equally the their memory and the suggested flow, and when they weight more their memory (and less the suggestions).
    \item $\initstate = (0.3, 0.5, 0.2)^\top$ and $\initstate = (0, 1, 0)^\top$: with all the other parameters fixed as $A = B = I_3, Q = Q_f = \begin{pmatrix}
         1 & 0 & 0 \\
         0 & 2 & 0 \\
         0 & 0 & 4
    \end{pmatrix},  R = \begin{pmatrix}
         0 & 0 & 0 \\
         0 & 0 & 0 \\
         0 & 0 & 0
    \end{pmatrix}, \gamma = 0.5$. This means we investigate the differences between the players' routing decision process when they start with putting partial flow on each edge, and when they start with putting all flow on the second edge.
    \item $Q = Q_f = I_3$ and $Q = Q_f = \begin{pmatrix}
         1 & 0 & 0 \\
         0 & 2 & 0 \\
         0 & 0 & 4
    \end{pmatrix}$: with all the other parameters fixed as $A = B = I_3, R =\begin{pmatrix}
         0 & 0 & 0 \\
         0 & 0 & 0 \\
         0 & 0 & 0
    \end{pmatrix}, \gamma = 0.5, \initstate = (0.3, 0.5, 0.2)^\top$. This means we investigate the differences between the players' routing decision process when the central decision maker penalizes equally the flow on the three edges, and when it gradually doubles the penalty on the edges.
\end{enumerate}

We also discuss the \textit{reliability of the solver} by showing one example, where the solver outputs the correct steady state solution only with large enough time horizon $T$. This example serves as an alert that in order to get the accurate numerical solution with the $\mpqpAlgo(\cdot)$, one needs to set a large time horizon $T$, at the price of longer computation time.


\subsection{Comparisons and Results}
For ease of reference, we list below the detailed choice of parameters in each figure.
In all experiments we assume no penalty on the controller (i.e., the central decision maker has the freedom to pose any suggested routing decision).\\

\begin{figure}[h]
\includegraphics[width=8cm]{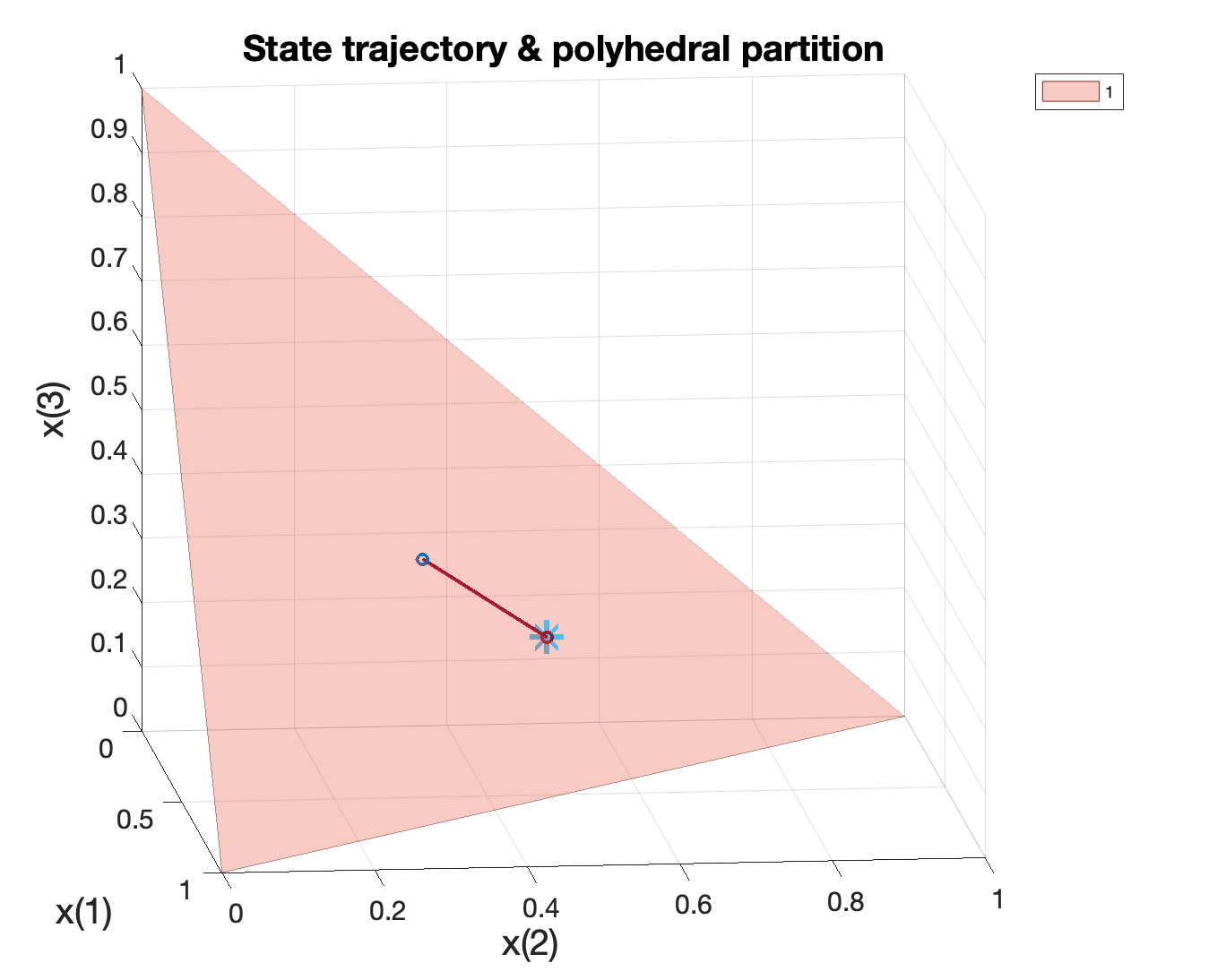}
\includegraphics[width=8cm]{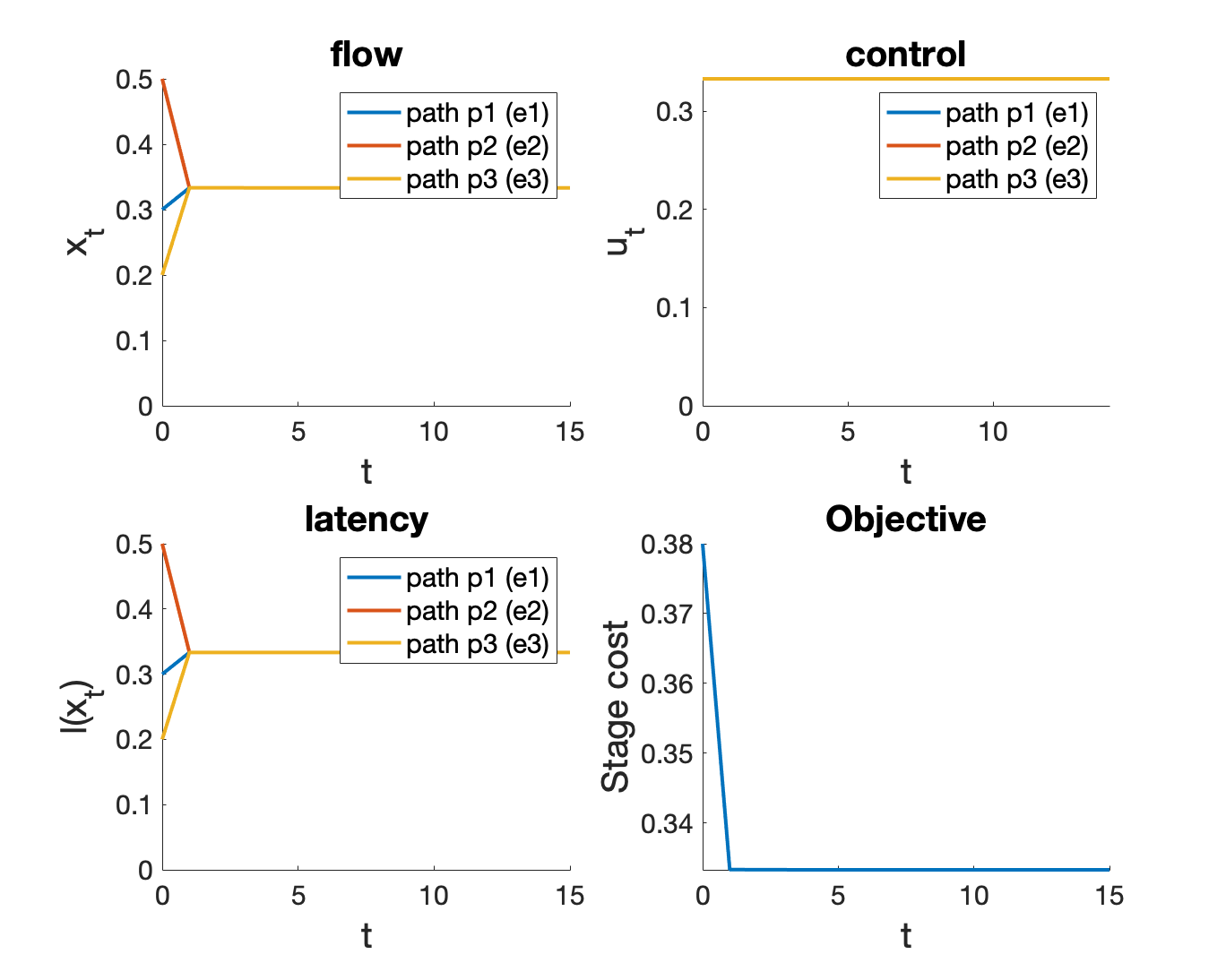}
\caption{Identity $Q, Q_f, B$, $A = (1/3), \gamma=0.5$. "*" denotes the initial state.}
\label{fig:A1:3_Q111_R0_gamma0.5}
\end{figure}

Figure \eqref{fig:A1:3_Q111_R0_gamma0.5}:
\begin{linenomath*}
\begin{align*}
 Q = Q_f = \begin{pmatrix}
         1 & 0 & 0 \\
         0 & 1 & 0 \\
         0 & 0 & 1
    \end{pmatrix},
    R = \begin{pmatrix}
         0 & 0 & 0 \\
         0 & 0 & 0 \\
         0 & 0 & 0
    \end{pmatrix}, 
A = &\begin{pmatrix}
         \frac{1}{3} & \frac{1}{3} & \frac{1}{3} \\
         \frac{1}{3} & \frac{1}{3} & \frac{1}{3} \\
         \frac{1}{3} & \frac{1}{3} & \frac{1}{3}
    \end{pmatrix},
    B = \begin{pmatrix}
         1 & 0 & 0 \\
         0 & 1 & 0 \\
         0 & 0 & 1
    \end{pmatrix},\\
    &T = 15,\gamma = 0.5, \initstate = \begin{pmatrix}
    0.3& 0.5&0.2\end{pmatrix}^\top
\end{align*}
\end{linenomath*}

\begin{figure}[h]
\includegraphics[width=8cm]{figures/A111_Q111_R0_gamma0.5_traj.png}
\includegraphics[width=8cm]{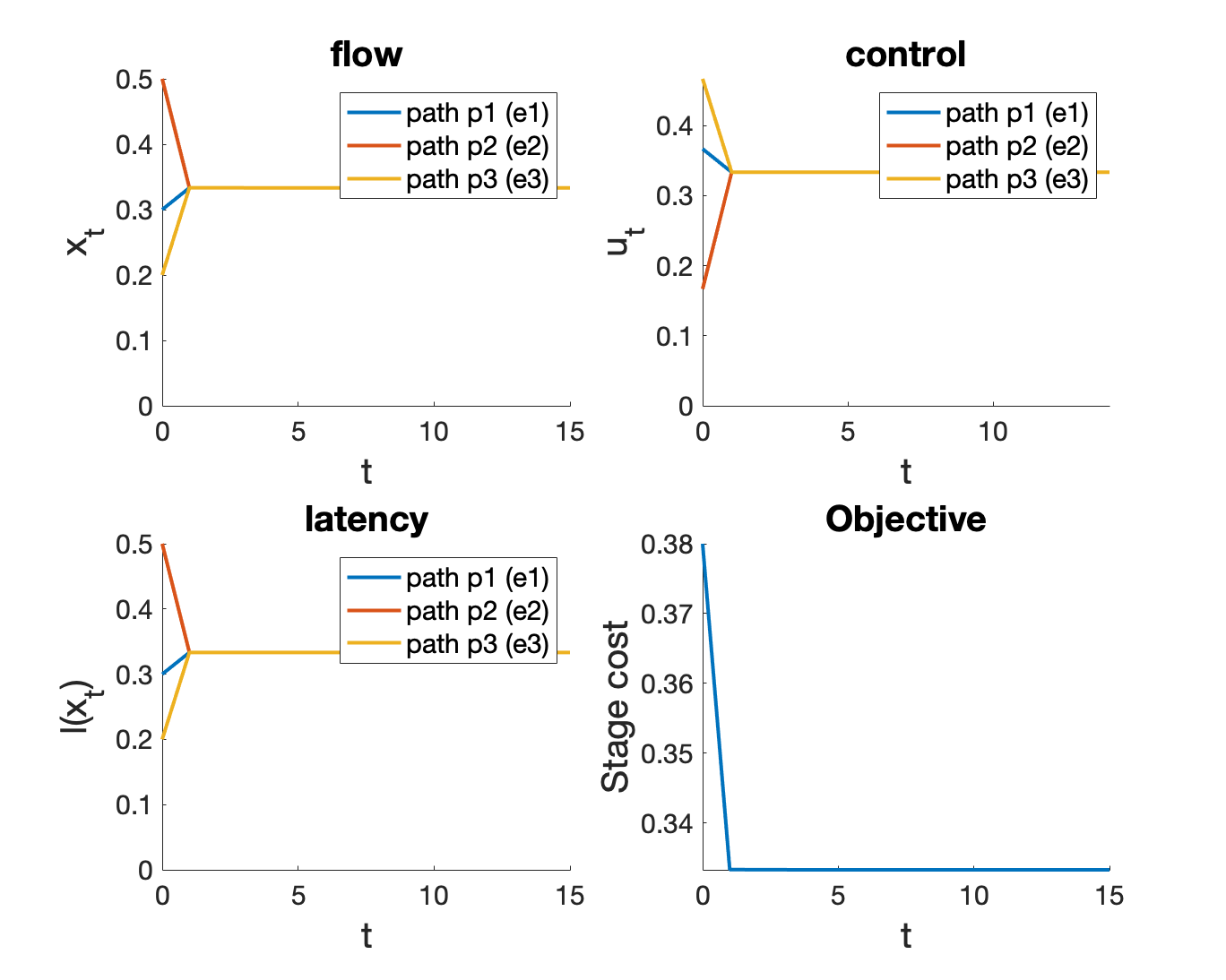}
\caption{Identity $Q, Q_f, A, B$, $\gamma=0.5$. "*" denotes the initial state.}
\label{fig:A111_Q111_R0_gamma0.5}
\end{figure}

Figure \eqref{fig:A111_Q111_R0_gamma0.5}:
\begin{linenomath*}
\begin{align*}
 Q = Q_f = \begin{pmatrix}
         1 & 0 & 0 \\
         0 & 1 & 0 \\
         0 & 0 & 1
    \end{pmatrix},
    R = \begin{pmatrix}
         0 & 0 & 0 \\
         0 & 0 & 0 \\
         0 & 0 & 0
    \end{pmatrix}, 
A = &\begin{pmatrix}
         1 & 0 & 0 \\
         0 & 1 & 0 \\
         0 & 0 & 1
    \end{pmatrix},
    B = \begin{pmatrix}
         1 & 0 & 0 \\
         0 & 1 & 0 \\
         0 & 0 & 1
    \end{pmatrix},\\
    &T = 15,\gamma = 0.5, \initstate = \begin{pmatrix}
    0.3& 0.5&0.2\end{pmatrix}^\top
\end{align*}
\end{linenomath*}

\begin{figure}[h]
\includegraphics[width=8cm]{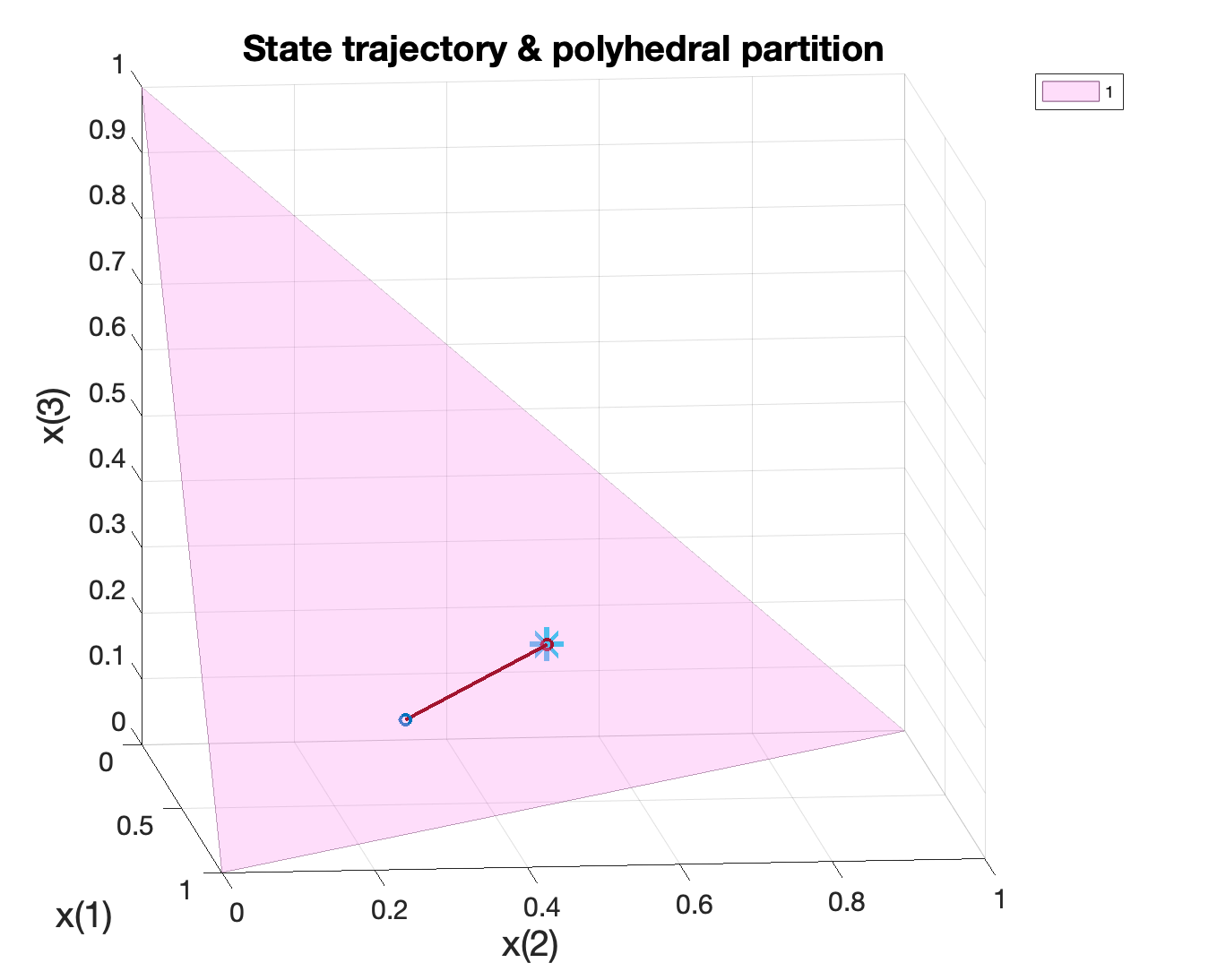}
\includegraphics[width=8cm]{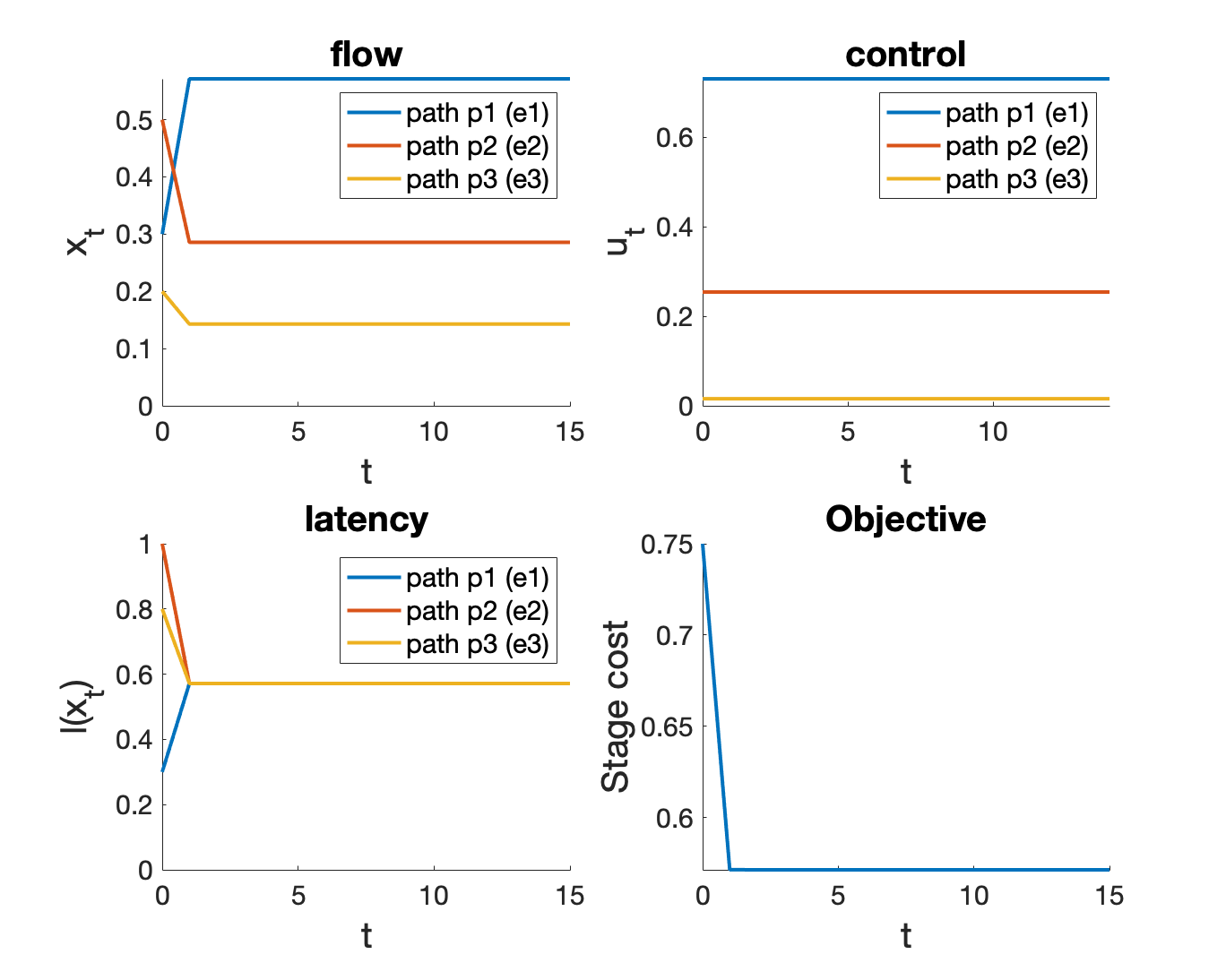}
\caption{Identity $B$, $A = (1/3), Q = Q_f = \mathrm{diag}(1,2,4),\gamma=0.5$. "*" denotes the initial state.}
\label{fig:A1:3_Q124_R0_gamma0.4}
\end{figure}

Figure \eqref{fig:A1:3_Q124_R0_gamma0.4}:
\begin{linenomath*}
\begin{align*}
 Q = Q_f = \begin{pmatrix}
         1 & 0 & 0 \\
         0 & 2 & 0 \\
         0 & 0 & 4
    \end{pmatrix},
    R = \begin{pmatrix}
         0 & 0 & 0 \\
         0 & 0 & 0 \\
         0 & 0 & 0
    \end{pmatrix}, 
A =& \begin{pmatrix}
         \frac{1}{3} & \frac{1}{3} & \frac{1}{3} \\
         \frac{1}{3} & \frac{1}{3} & \frac{1}{3} \\
         \frac{1}{3} & \frac{1}{3} & \frac{1}{3}
    \end{pmatrix},
    B = \begin{pmatrix}
         1 & 0 & 0 \\
         0 & 1 & 0 \\
         0 & 0 & 1
    \end{pmatrix},\\
    &T = 15,\gamma = 0.5, \initstate = \begin{pmatrix}
    0.3& 0.5&0.2\end{pmatrix}^\top
\end{align*}
\end{linenomath*}

\begin{figure}[h]
\includegraphics[width=8cm]{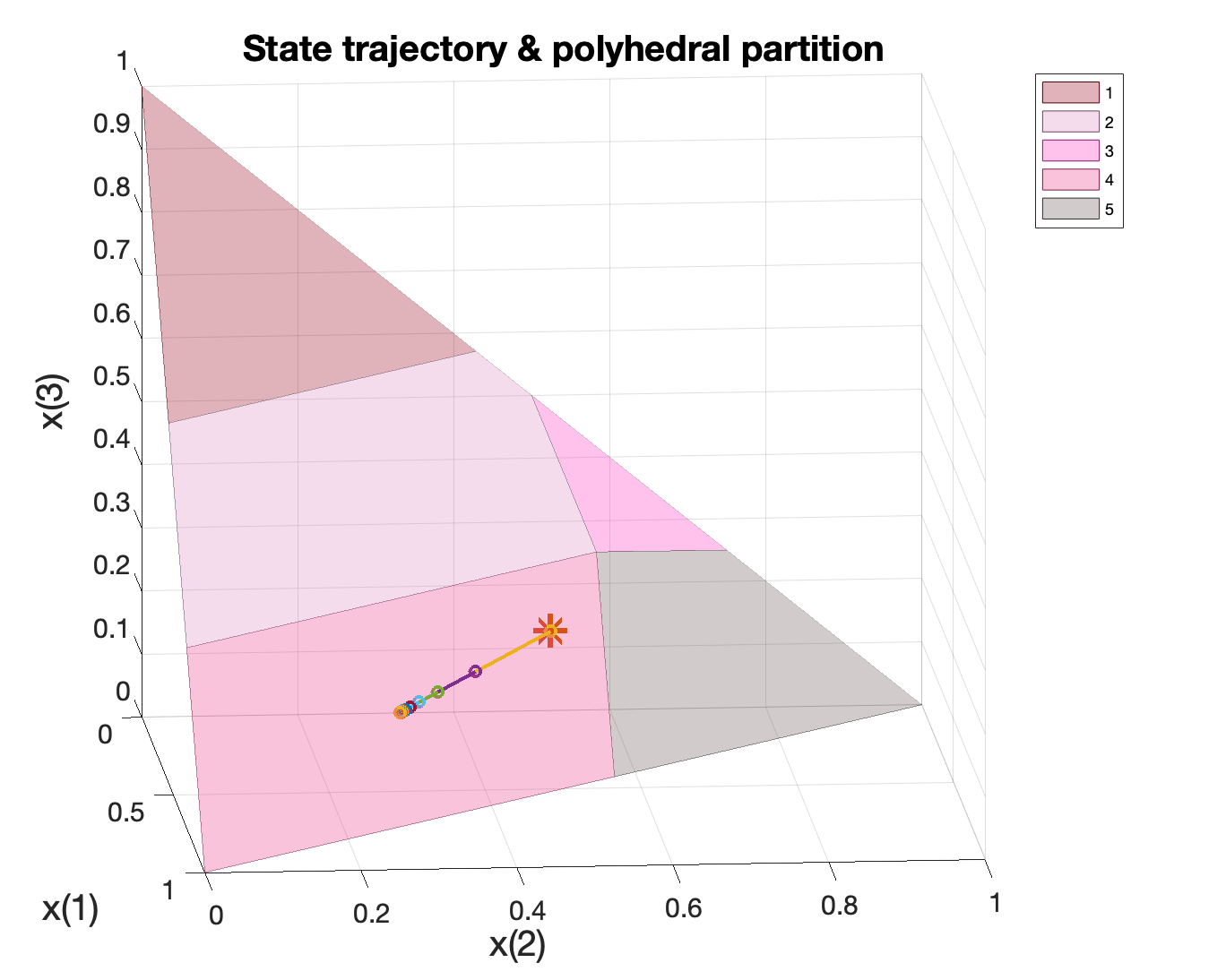}
\includegraphics[width=8cm]{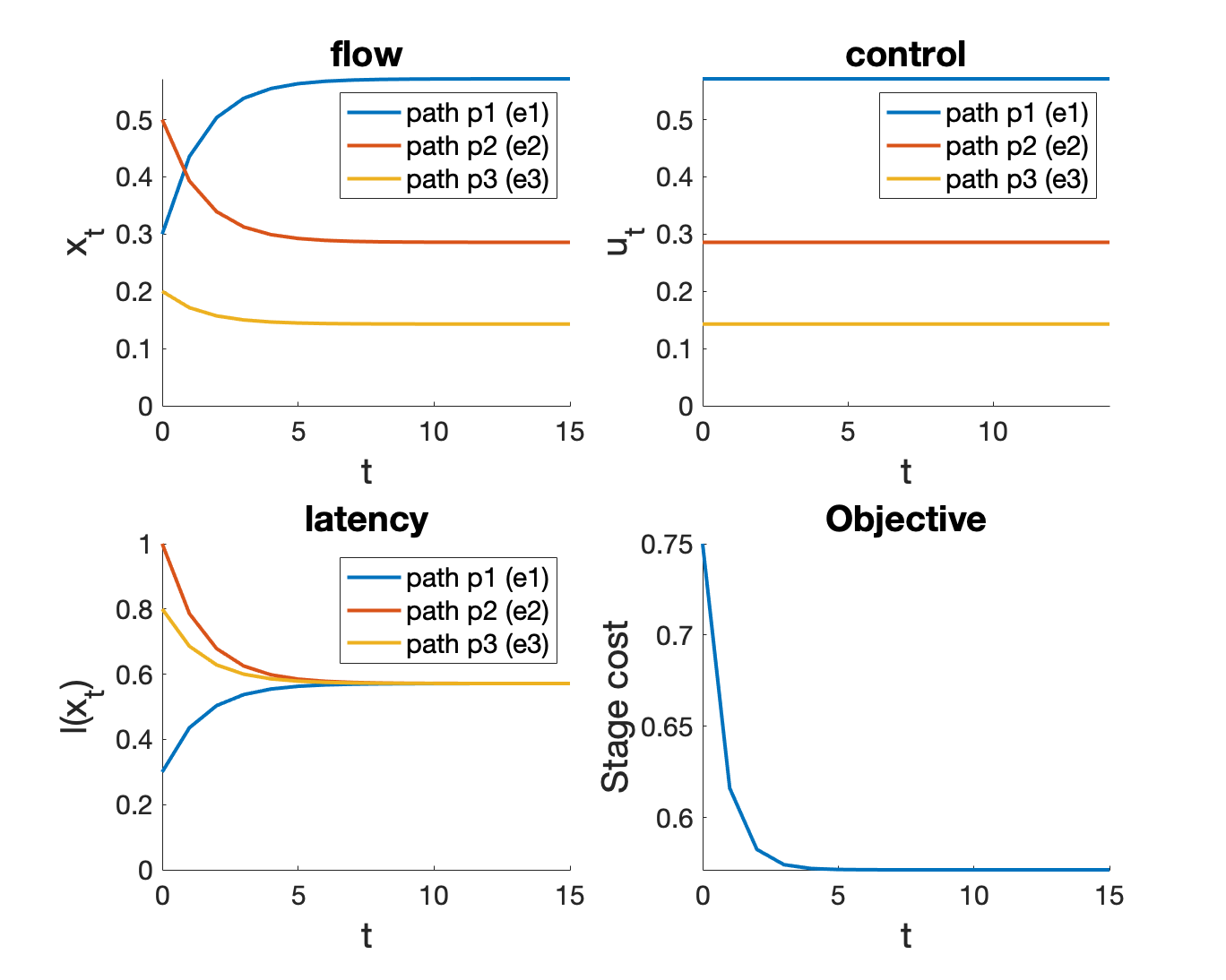}
\caption{Identity $A,B$, $Q = Q_f = \mathrm{diag}(1,2,4),\gamma=0.5$. "*" denotes the initial state.}
\label{fig:A111_Q124_R0_gamma0.5}
\end{figure}

Figure \eqref{fig:A111_Q124_R0_gamma0.5}:
\begin{linenomath*}
\begin{align*}
 Q = Q_f = \begin{pmatrix}
         1 & 0 & 0 \\
         0 & 2 & 0 \\
         0 & 0 & 4
    \end{pmatrix},
    R = \begin{pmatrix}
         0 & 0 & 0 \\
         0 & 0 & 0 \\
         0 & 0 & 0
    \end{pmatrix}, 
A =& \begin{pmatrix}
         1 & 0 & 0 \\
         0 & 1 & 0 \\
         0 & 0 & 1
    \end{pmatrix},
    B = \begin{pmatrix}
         1 & 0 & 0 \\
         0 & 1 & 0 \\
         0 & 0 & 1
    \end{pmatrix},\\
    &T = 15,\gamma = 0.5, \initstate = \begin{pmatrix}
    0.3& 0.5&0.2\end{pmatrix}^\top
\end{align*}
\end{linenomath*}

\begin{figure}[h]
\includegraphics[width=8cm]{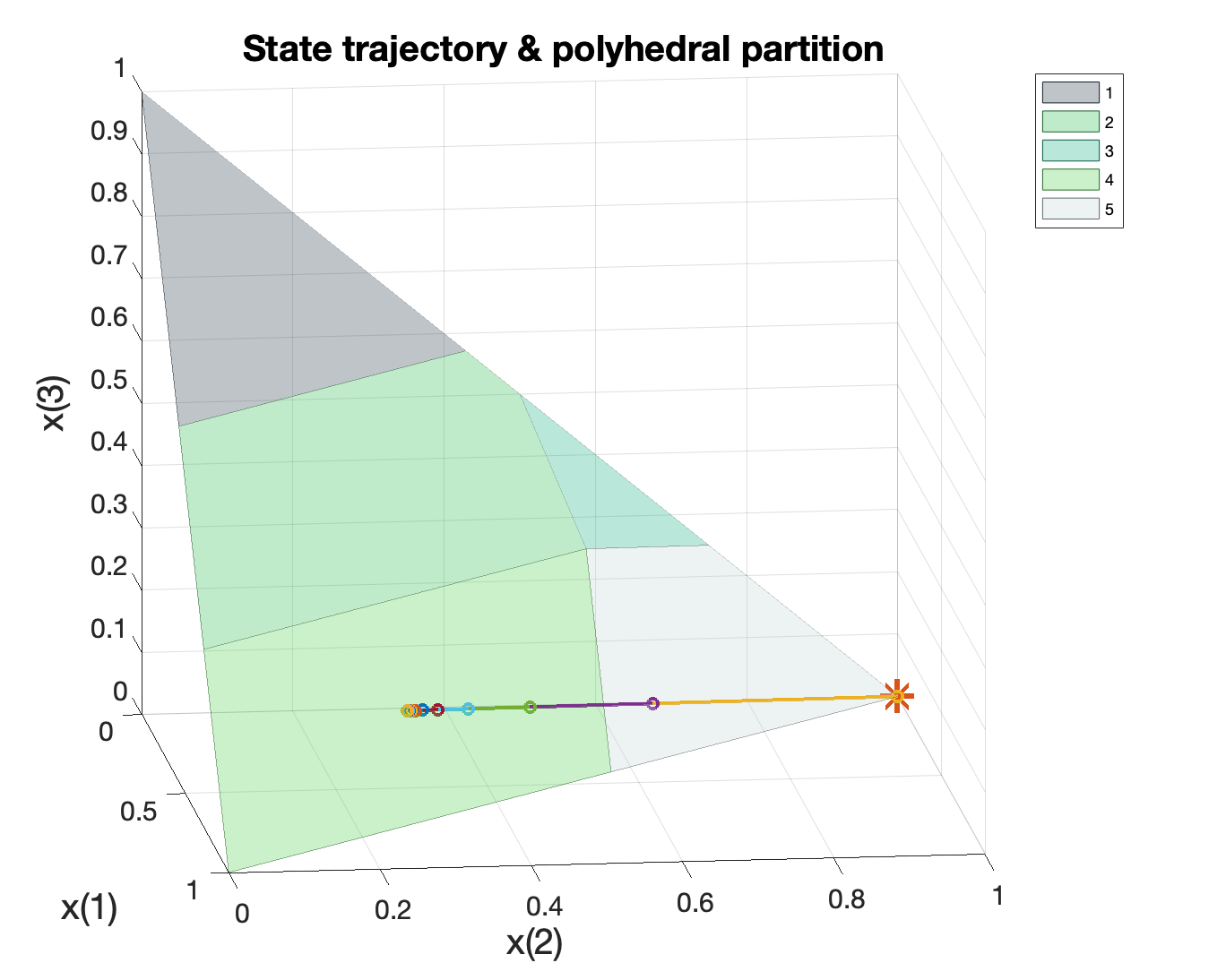}
\includegraphics[width=8cm]{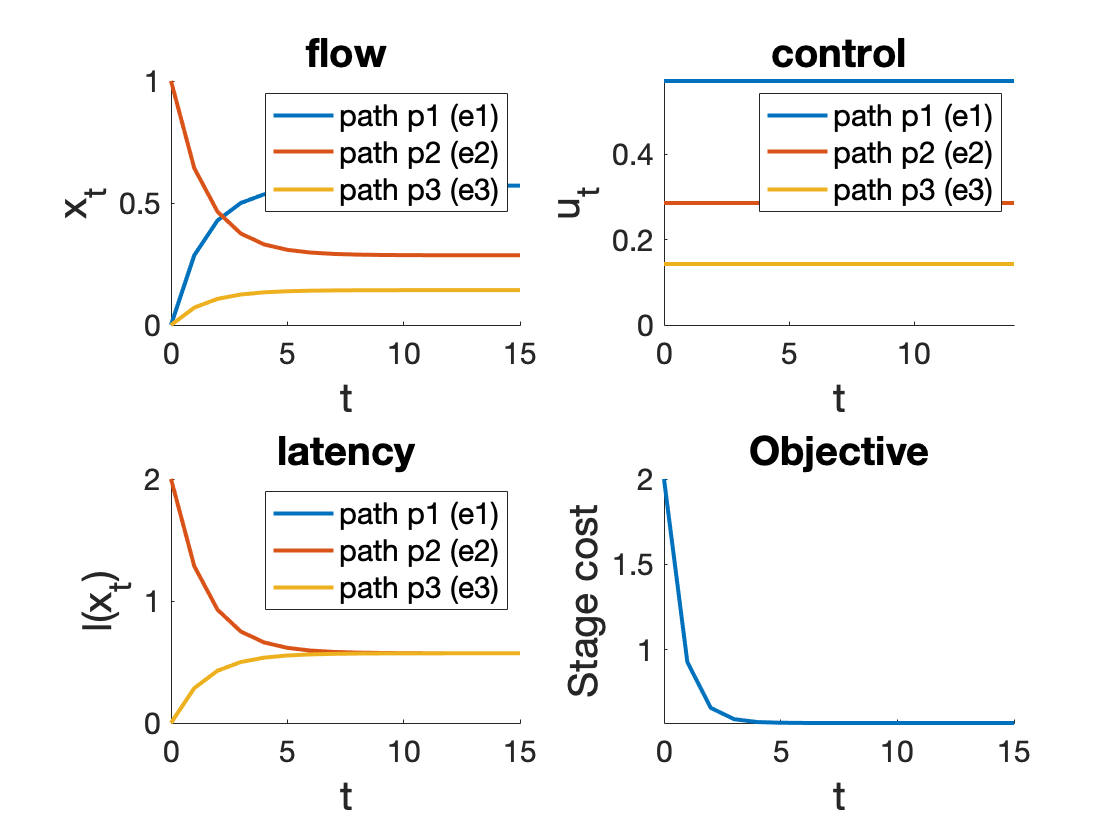}
\caption{Identity $A,B$, $Q=Q_f=\mathrm{diag}(1,2,4),\gamma=0.5$, initial $(0,1,0)$. "*" denotes the initial state.}
\label{fig:A111_Q124_R0_gamma0.5_initial}
\end{figure}

Figure \eqref{fig:A111_Q124_R0_gamma0.5_initial}:
\begin{linenomath*}
\begin{align*}
 Q = Q_f = \begin{pmatrix}
         1 & 0 & 0 \\
         0 & 2 & 0 \\
         0 & 0 & 4
    \end{pmatrix},
    R = \begin{pmatrix}
         0 & 0 & 0 \\
         0 & 0 & 0 \\
         0 & 0 & 0
    \end{pmatrix}, 
A =& \begin{pmatrix}
         1 & 0 & 0 \\
         0 & 1 & 0 \\
         0 & 0 & 1
    \end{pmatrix},
    B = \begin{pmatrix}
         1 & 0 & 0 \\
         0 & 1 & 0 \\
         0 & 0 & 1
    \end{pmatrix},\\
    &T = 15,\gamma = 0.5, \initstate = \begin{pmatrix}
    0&1&0\end{pmatrix}^\top
\end{align*}
\end{linenomath*}

\begin{figure}[h]
\includegraphics[width=8cm]{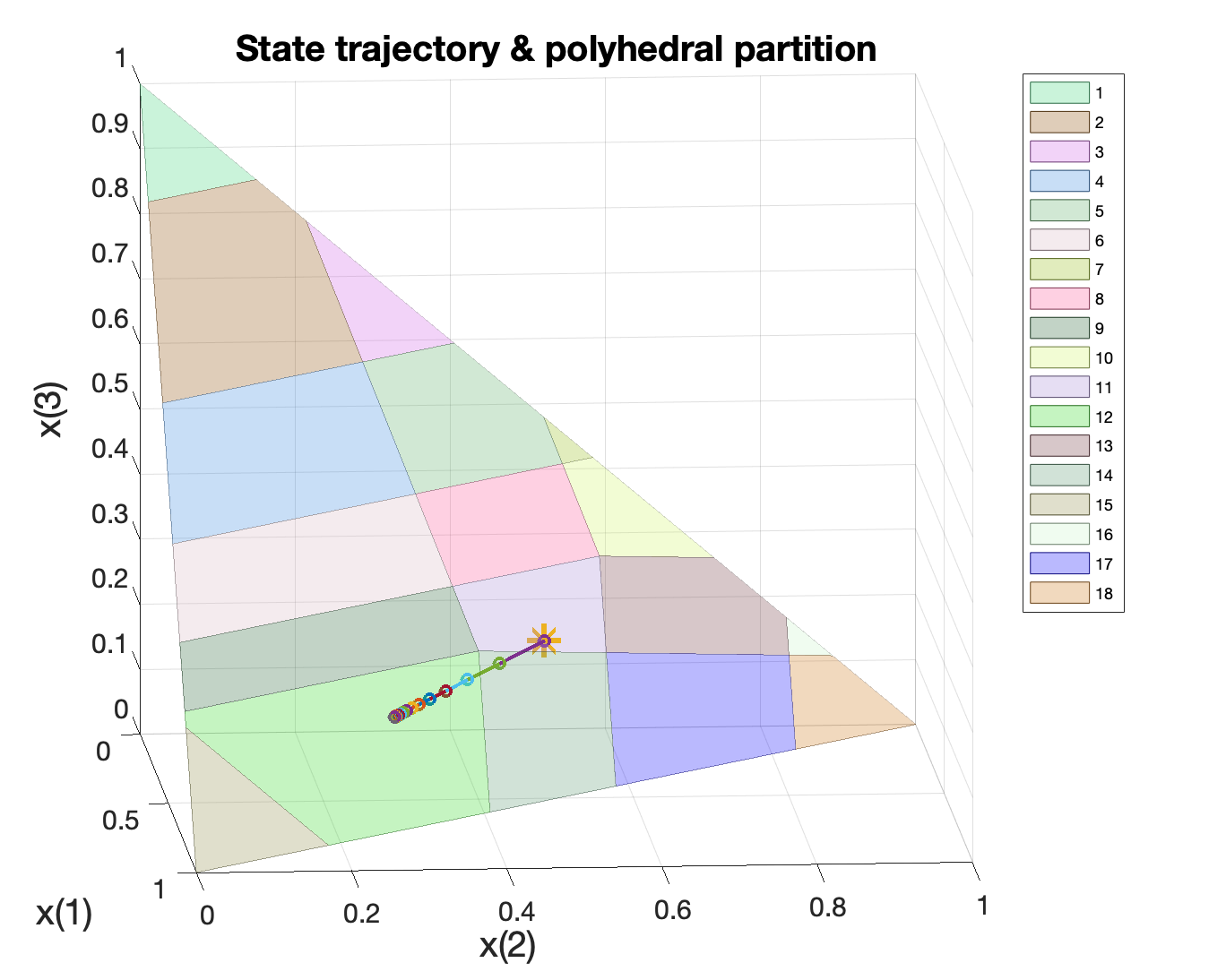}
\includegraphics[width=8cm]{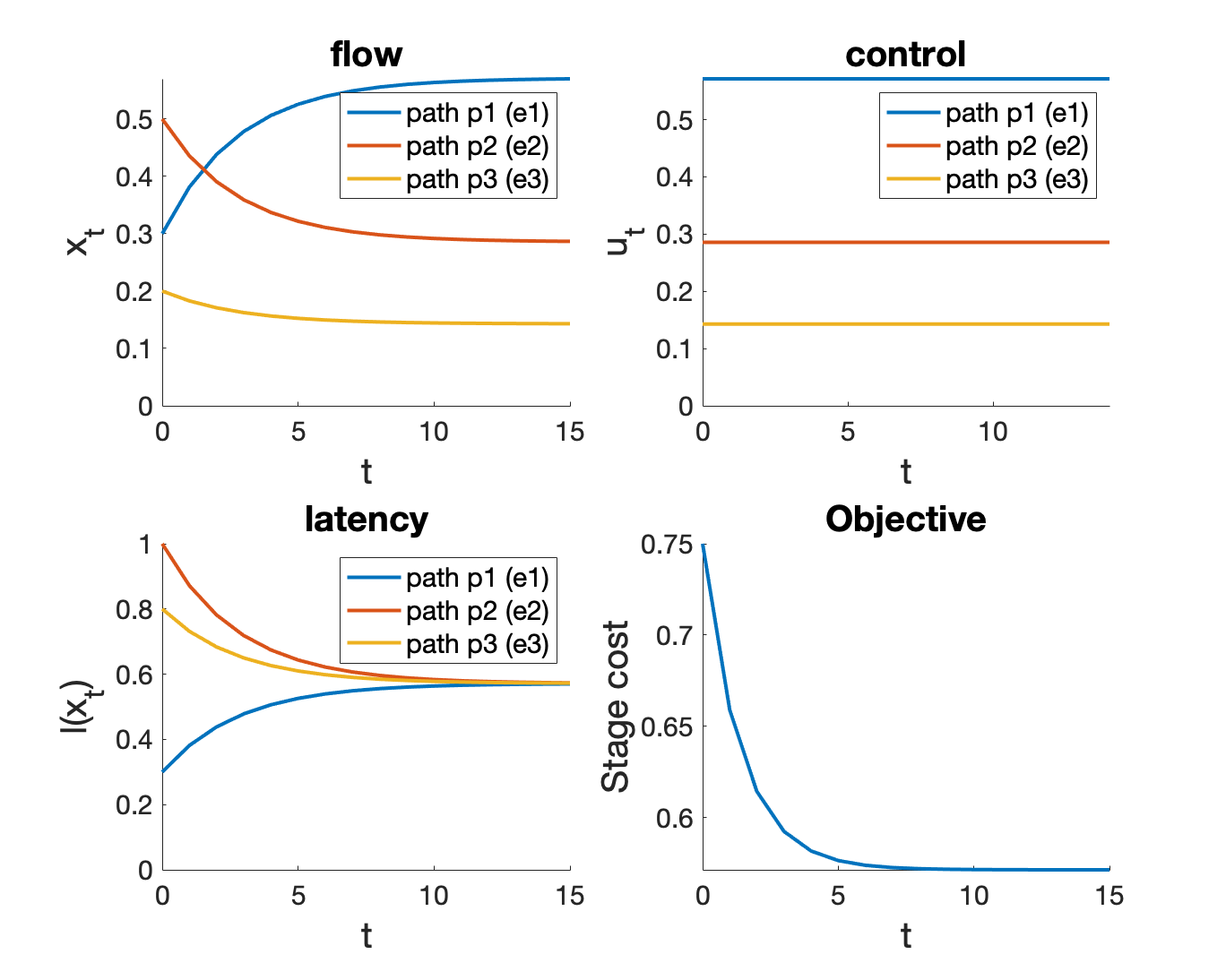}
\caption{Identity $A,B$, $Q = Q_f = \mathrm{diag}(1,2,4),\gamma=0.7$. "*" denotes the initial state.}
\label{fig:A111_Q124_R0_gamma0.7}
\end{figure}

Figure \eqref{fig:A111_Q124_R0_gamma0.7}:
\begin{linenomath*}
\begin{align*}
 Q = Q_f = \begin{pmatrix}
         1 & 0 & 0 \\
         0 & 2 & 0 \\
         0 & 0 & 4
    \end{pmatrix},
    R = \begin{pmatrix}
         0 & 0 & 0 \\
         0 & 0 & 0 \\
         0 & 0 & 0
    \end{pmatrix}, 
A =& \begin{pmatrix}
         1 & 0 & 0 \\
         0 & 1 & 0 \\
         0 & 0 & 1
    \end{pmatrix},
    B = \begin{pmatrix}
         1 & 0 & 0 \\
         0 & 1 & 0 \\
         0 & 0 & 1
    \end{pmatrix},\\
    &T = 15,\gamma = 0.7, \initstate = \begin{pmatrix}
    0.3&0.2&0.5\end{pmatrix}^\top
\end{align*}
\end{linenomath*}

\textbf{Discussion:}
In each figure, we present on the left the probability simplex in $\mathbb{R}^3$, the critical regions in it (recall that within each partition the optimal controller takes different affine expression), and the convergence trajectory of the state $\stateiter$ (with the initial state denoted by $*$). On the right we present (1) the value of the flow allocation $\stateiter$ on each edge at each time step $\timesym$, (2) the value of the optimal control (suggested routing decision) $\controliter$ on each edge at each time step $\timesym$, (3) the incurred latency $\lossiter$ on each edge at each time step $\timesym$, and (4) the incurred stage cost $\stateiter^\top Q \stateiter + \controliter^\top R \controliter, t = 0, \dots, T-1$ and the terminal cost $\state_T^\top Q_f \state_T$.
\begin{enumerate}
    \item \textbf{Varying $A$ matrices in the dynamics:} To see the differences in changing $A$, we compare figure \eqref{fig:A1:3_Q111_R0_gamma0.5} and figure \eqref{fig:A111_Q111_R0_gamma0.5}.
    \begin{itemize}
        \item It can be observed that in both scenarios, the controlled system converges to its steady state $\state\optimal = (1/3, 1/3, 1/3)^\top$ in one step. The numerical results suggest that with the choice of cost matrices being $Q = Q_f = I_3, R = 0$ and the specific setup of other parameters in the numerics, the steady state $\state\optimal$ coincides with the Nash equilibrium in the parallel network with latency functions $l_1(y) = l_2(y) = l_3(y) = y$. 
        \item Specifically, the fact that the state $\mathbf{x}_t$ converges in one step to the Nash equilibrium with $A = \begin{pmatrix}
         \frac{1}{3} & \frac{1}{3} & \frac{1}{3} \\
         \frac{1}{3} & \frac{1}{3} & \frac{1}{3} \\
         \frac{1}{3} & \frac{1}{3} & \frac{1}{3}
    \end{pmatrix}$ agrees with our analysis in property (\ref{1/N_A}).
        \item The polyhedral partitions of the feasible state space $\flowset$ are different. In the case of $A = \begin{pmatrix}
         \frac{1}{3} & \frac{1}{3} & \frac{1}{3} \\
         \frac{1}{3} & \frac{1}{3} & \frac{1}{3} \\
         \frac{1}{3} & \frac{1}{3} & \frac{1}{3}
    \end{pmatrix}$ there is no partition of the space, meaning that the optimal controller takes the same affine function over $\flowset$. In the case of $A = I_3$, the feasible state space is partitioned into four polyhedral regions, within each taking different expression of affine optimal controllers.
    \end{itemize}
    \item \textbf{Varying $\gamma$ in the dynamics:} To see the differences in changing $\gamma$, we compare figure \eqref{fig:A111_Q124_R0_gamma0.5} and figure \eqref{fig:A111_Q124_R0_gamma0.7}.
    \begin{itemize}
        \item Both controlled systems converge to the steady state $\state\optimal = (4/7, 2/7, 1/7)^\top$. The numerical results suggest that with the choice of cost matrices being $Q = Q_f = \begin{pmatrix}
         1 & 0 & 0 \\
         0 & 2 & 0 \\
         0 & 0 & 4
    \end{pmatrix}, R = 0$ and the specific setup of other parameters in the numerics, the steady state $\state\optimal$ coincides with the Nash equilibrium in the parallel network with latency functions $l_1(y) = y, l_2(y) = 2y, l_3(y) = 4y$. 
    \item The polyhedral partitions of the feasible state space $\Delta$ are different in the two cases. We observe that with larger $\gamma = 0.7$ (i.e. the players update their flow allocation considering \textbf{more} their memory than the central decision maker's suggested routing decision), the feasible state space has way more polyhedral partitions. This intuitively means that: if the players weight less the suggestions from the central decision maker, in order for the central decision maker to drive the traffic to a target flow allocation, it will have to design more complex routing strategy based on the actual flow on the network.
    \item The convergence rate with larger $\gamma$ is slower: we see that with $\gamma = 0.5$, the state stabilizes with around 5 time steps, while with $\gamma = 0.7$, the state stabilizes with around 10 time steps. This can be intuitively interpreted as: if the players weight less the suggestions from the central decision maker, the steady state of the system is less controllable and therefore it takes longer for the central decision maker to stabilize the system.
    \end{itemize}
    \item \textbf{Varying initial states $\initstate$:} To see the differences in changing $\initstate$, we compare figure \eqref{fig:A111_Q124_R0_gamma0.5} and figure \eqref{fig:A111_Q124_R0_gamma0.5_initial}.
    \begin{itemize}
        \item Both controlled systems converge to the steady state $\state\optimal = (4/7, 2/7, 1/7)^\top$, which again coincides with the Nash equilibrium in the parallel network with latency functions $l_1(y) = y, l_2(y) = 2y, l_3(y) = 4y$. 
        \item We would like to point out that with ${x}_0 = [0,1,0]^\top$, the state of the system starts further away from its steady state, and starts within a different polyhedron from where the steady state is located. However with the specific control dynamics in the numerical experiment, the central decision maker still manages to steer the system to the target flow allocation.
    \end{itemize}
    \item \textbf{Varying cost functions:} To see the differences in changing the cost matrices, we compare figure \eqref{fig:A111_Q111_R0_gamma0.5} and figure \eqref{fig:A111_Q124_R0_gamma0.5}.
    \begin{itemize}
        \item Both controlled systems converge, but to the different steady states. More specifically, with $Q = Q_f = I_3$, the steady state is $\state\optimal = (1/3, 1/3, 1/3)^\top$, which coincides with the Nash equilibrium in the parallel network with latency functions $l_1(y) = l_2(y) = l_3(y) = y$; with $Q = Q_f = \begin{pmatrix}
         1 & 0 & 0 \\
         0 & 2 & 0 \\
         0 & 0 & 4
    \end{pmatrix}$, the steady state is $\state\optimal = (4/7, 2/7, 1/7)^\top$, which coincides with the Nash equilibrium in the parallel network with latency functions $l_1(y) = y, l_2(y) = 2y, l_3(y) = 4y$.
    \item This justifies the capability of the central decision maker to design different objective functions and steer the controlled system to different target flow allocations.
    \end{itemize}
\end{enumerate}

\textbf{Reliability of the solver:} with the following setting, we show that it is important to set the time horizon $T$ to be large enough so that the algorithm gives the steady state solution. \\

\begin{figure}[h]
\includegraphics[width=8cm]{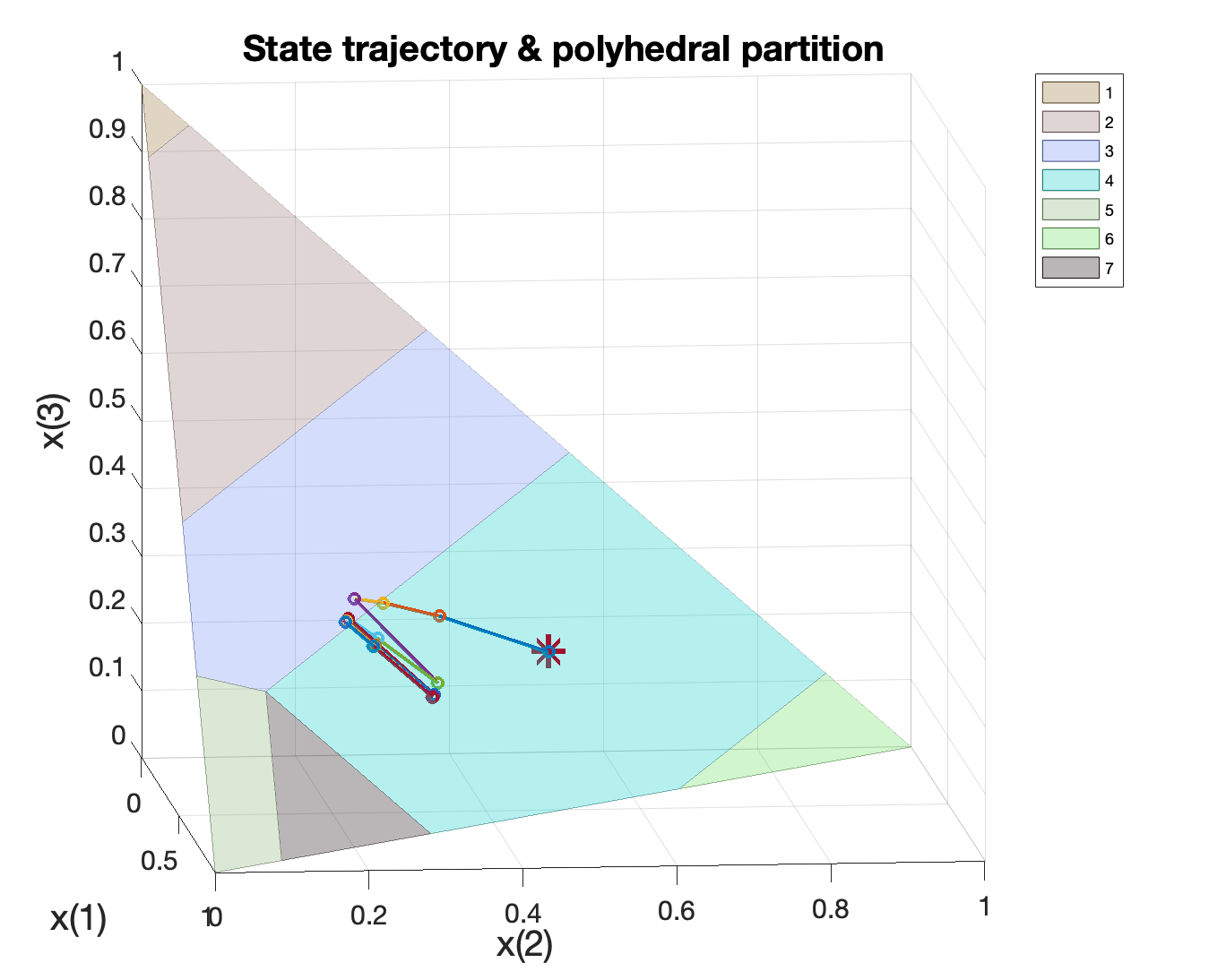}
\includegraphics[width=8cm]{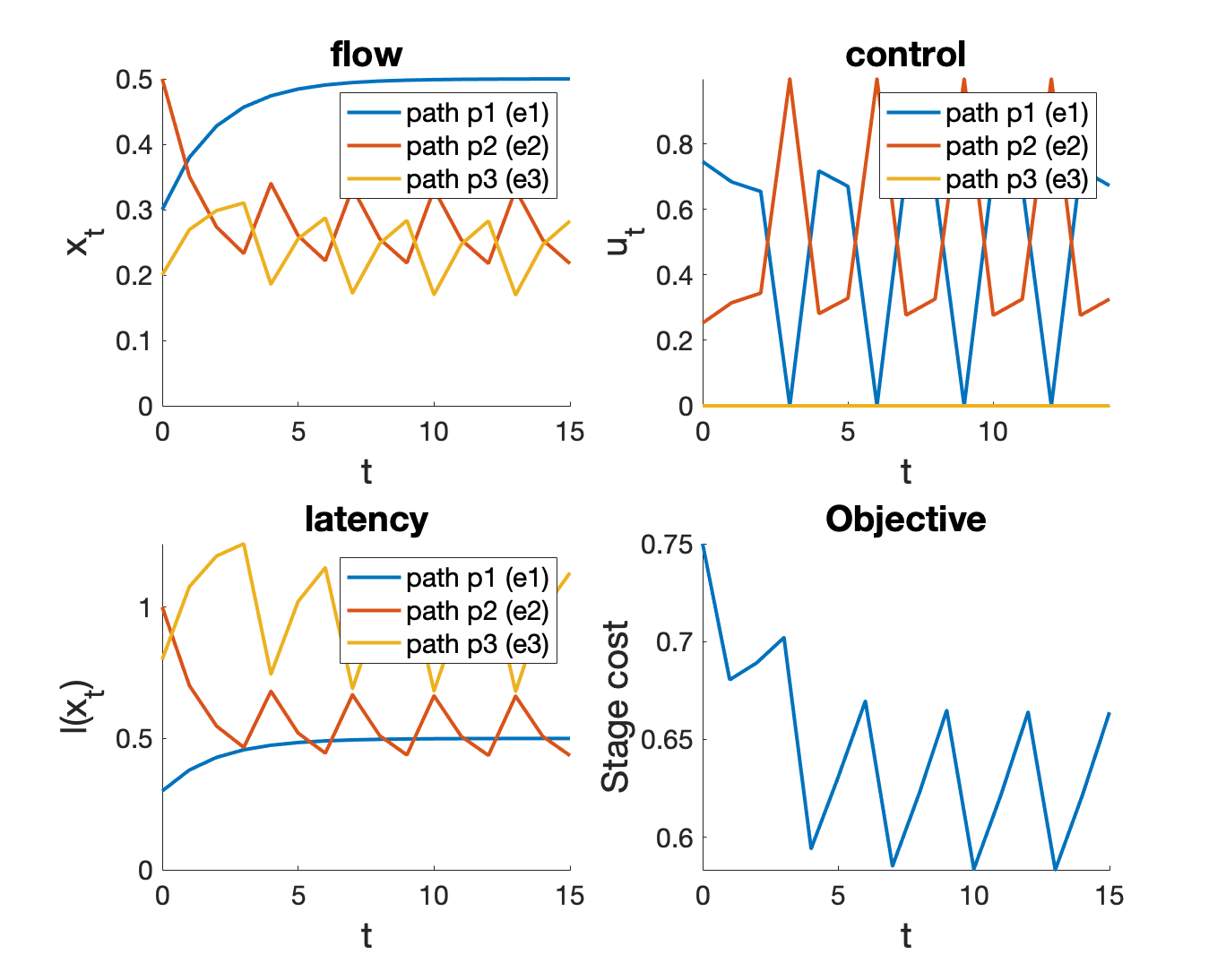}
\caption{Identity $A$, $B =\begin{pmatrix}
         1/2 & 1/2 & 0 \\
         0 & 1/2 & 1/2 \\
         1/2 & 0 & 1/2
    \end{pmatrix}$, $Q=Q_f=\mathrm{diag}(1,2,4),\gamma=0.6$, $T = 15$. "*" denotes the initial state.}
\label{fig:Non_coverge}
\end{figure}

\begin{figure}[h]
\includegraphics[width=8cm]{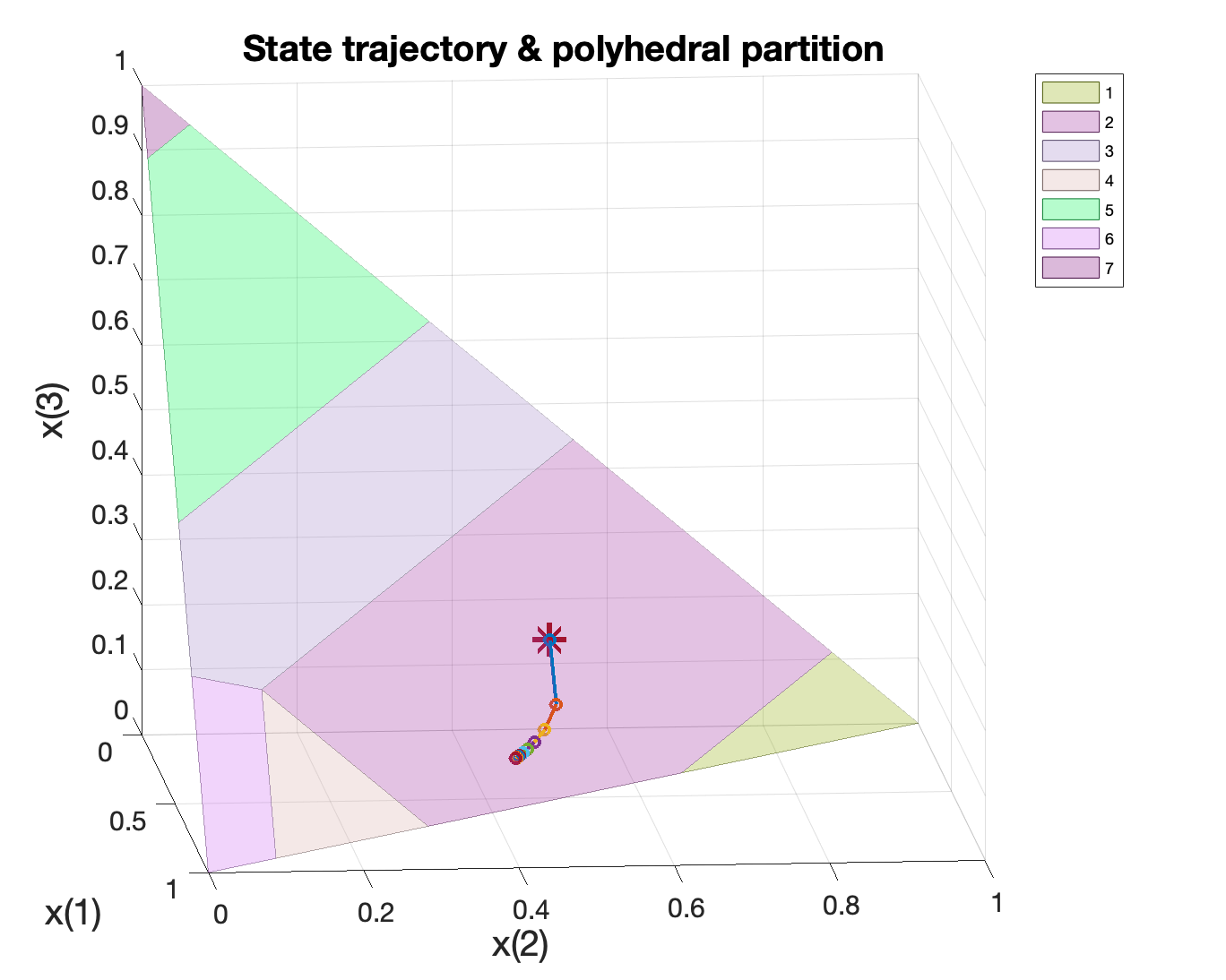}
\includegraphics[width=8cm]{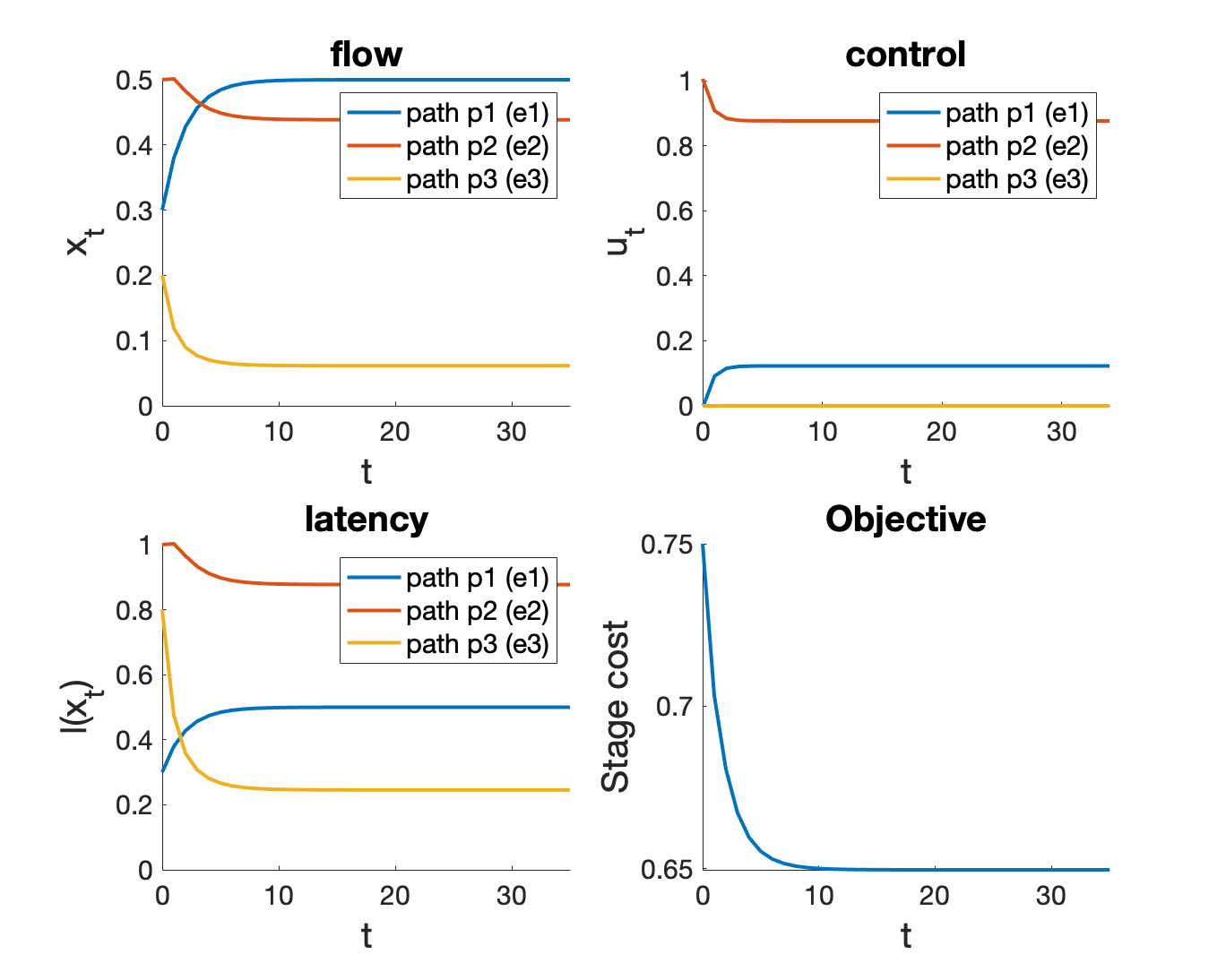}
\caption{Identity $A$, $B =\begin{pmatrix}
         1/2 & 1/2 & 0 \\
         0 & 1/2 & 1/2 \\
         1/2 & 0 & 1/2
    \end{pmatrix}$, $Q=Q_f=\mathrm{diag}(1,2,4),\gamma=0.6$, $T = 35$. "*" denotes the initial state.}
\label{fig:coverge}
\end{figure}

Let
\begin{linenomath*}
\begin{align*}
 Q = Q_f = \begin{pmatrix}
         1 & 0 & 0 \\
         0 & 2 & 0 \\
         0 & 0 & 4
    \end{pmatrix},
    R = \begin{pmatrix}
         0 & 0 & 0 \\
         0 & 0 & 0 \\
         0 & 0 & 0
    \end{pmatrix},
A =& \begin{pmatrix}
         1 & 0 & 0 \\
         0 & 1 & 0 \\
         0 & 0 & 1
    \end{pmatrix},
    B = \begin{pmatrix}
         \frac{1}{2} & \frac{1}{2} & 0 \\
         0 & \frac{1}{2} & \frac{1}{2} \\
         \frac{1}{2} & 0 & \frac{1}{2}
    \end{pmatrix},
    \\
    &\gamma = 0.6, \initstate = \begin{pmatrix}
    0.3&0.5&0.2\end{pmatrix}^\top
\end{align*}
\end{linenomath*}
If we set the time horizon $T = 15$, we observe in figure \eqref{fig:Non_coverge} that the control system can be oscillating and non-stable, and does not converge to a steady state. However if we set a larger $T = 35$, we see in figure \eqref{fig:coverge} that the controlled system is stabilizing -- the trajectory of the state converges to a steady state.
\\

\section{Conclusion}
This article transformed the dynamic traffic assignment problem into a control theoretic problem, under a simple scenario using a parallel network with affine latency functions on each edge. We started with the routing problem under the repeated game framework, proposed a game design scheme through LQR within the control theoretic framework, and then leveraged control techniques to analyze convergence/stability of the system with special examples. Within our system analysis, we discovered that the design of $A, B$ in the equation of the dynamics could be chosen to be \textit{left stochastic matrices} so that the conservation of flow is guaranteed. Additionally, the suggested routing allocation $\control$ has to be chosen such that the resulting vector is stochastic, disqualifying certain control technique that can be used in this system (i.e. the algebraic Ricatti equation) and requesting the use of others. Using this fact, we designed an algorithmic solution for the optimal control using a multiparametric quadratic programming approach (explicit MPC), generating explicit strategies and geometric plots to visualize the optimal solutions the central decision maker should take at the feasible states of the game. We discussed the impact of each system parameter on the solution and explained that various structures in $A$ will result in different solution requirements for the central decision maker, while the cost function changes the system's steady state, and $\gamma, \initstate$ influence the rate of convergence. Future work includes studying the region of attraction of the LQR routing system via reachability analysis, extending the present work to more complex network structures (beyond the parallel network), generalizing the analysis of the convergence behavior of the LQR routing system beyond the special cases considered in the present work, and leveraging the Reinforcement Learning techniques to approach general routing problems.

\bibliographystyle{trb}
\bibliography{trb_template}
\end{document}